\ifx\shlhetal\undefinedcontrolsequence\let\shlhetal\relax\fi
\documentclass[11pt]{article}

\usepackage{amssymb}
\usepackage{amsfonts}
\usepackage{amsbsy}

\renewcommand{\lg}{{\ell g}}

\newcommand{\defeq}{\stackrel{\rm def}{=}}
\newcommand{\Op}{{\rm Op}}
\newcommand{\id}{{\rm id}}
\newcommand{\rrk}{{\rm rk}}
\newcommand{\otp}{{\rm otp}}
\newcommand{\Rang}{{\rm Rang}}

\newcommand{\bT}{{\bf T}}
\newcommand{\cf}{{\rm cf}}
\newcommand{\dcf}{{\rm dcf}}
\newcommand{\can}{{\rm canonical}}
\newcommand{\Lim}{{\rm Lim}}
\newcommand{\cF}{{\cal F}}
\newcommand{\cK}{{\cal K}}

\newcommand{\tp}{{\rm tp}}

\newcommand{\qed}{\hspace{0.2in}\vrule width 6pt height 6pt depth 0pt
\vspace{0.1in}}

\newcommand{\boldk}{{\boldsymbol{\kappa}}}
\newcommand{\prenice}{\preceq_{\rm nice}}
\newcommand{\nice}{{\rm nice}}
\newcommand{\imply}{\Rightarrow}

   \def\nofork^#1_#2{\putforkinmargin
        \unionstick^{\textstyle #1}_{\textstyle #2}}

\newbox\noforkbox \newdimen\forklinewidth
\setbox0\hbox{$\textstyle\bigcup$}
\setbox1\hbox to \wd0{\hfil\vrule width \forklinewidth depth \dp0
                        height \ht0 \hfil}
\setbox\noforkbox\hbox{\box1\box0\relax}
\def\unionstick{\mathop{\copy\noforkbox}\limits}
\def\nonfork#1#2_#3{#1\unionstick_{\textstyle #3}#2}
\def\nonforkin#1#2_#3^#4{#1\unionstick_{\textstyle #3}^{\textstyle #4}#2}

\def\doesforkin_#1^#2{\biguplus\limits_{\textstyle #1}^{\textstyle #2}}

\newcommand{\Proof}{{\sc Proof} \hspace{0.2in}}
\newtheorem{theorem}{Theorem}[section]
\newtheorem{data}{Data}[theorem]
\newtheorem{ftry}[data]{First Try}
\newtheorem{context}[theorem]{Context}
\newtheorem{fact}[theorem]{Fact}
\newtheorem{moddata}{Modified Data}[theorem]
\newtheorem{sclaim}[moddata]{Claim}
\newtheorem{sfact}[moddata]{Fact}
\newtheorem{improvement}[moddata]{Improvement}
\newtheorem{sremark}{Remark}[theorem]
\newtheorem{snotation}[sremark]{Notation}
\newtheorem{ssremark}[sremark]{Remark}
\newtheorem{ssfact}[sremark]{Fact}

\newtheorem{construction}[theorem]{Construction}
\newtheorem{claim}[theorem]{Claim}
\newtheorem{lemma}[theorem]{Lemma}
\newtheorem{discussion}[theorem]{Discussion}
\newtheorem{hypothesis}[theorem]{Hypothesis}

\newtheorem{corollary}[theorem]{Corollary}
\newtheorem{conclusion}[theorem]{Conclusion}
\newtheorem{intermediatecorollary}[theorem]{Intermediate Corollary}

\newtheorem{remark}[theorem]{Remark}

\newtheorem{definition}[theorem]{Definition}

\newtheorem{example}[theorem]{Example}

\date\today

\title{Categoricity of Theories in $L_{\kappa^* \omega}$, when \\
$\kappa^*$ is a measurable cardinal. Part II}

\author{Saharon Shelah\thanks{Research supported by the United
States-Israel Binational Science Foundation. 
Done 6-7/88. Publication number 472}\\
Institute of Mathematics\\
The Hebrew University of Jerusalem\\
91904 Jerusalem, Israel\\
and\\
Department of Mathematics\\
Rutgers University\\
New Brunswick, NJ 08854, USA
}

\begin{document}
\baselineskip13.14 truept

\maketitle

\begin{abstract}
We continue the work of \cite{KlSh:362} and prove that for $\lambda$
successor, a $\lambda$-categorical
theory $\bT$ in $L_{\kappa^*,\omega}$ is $\mu$-categorical for every $\mu$,
$\mu\leq\lambda$ which is above the $(2^{LS(\bT)})^+$-beth cardinal.
\end{abstract}
\vfill
\eject

\section{Introduction}
We deal here with the categoricity spectrum of theory $\bT$ in logic:
$\bT\subseteq L_{\kappa^*, \omega}$ with $\kappa^*$ measurable. Makkai
Shelah \cite{MaSh:285} have dealt with the case $\kappa^*$ a compact
cardinal. So $\kappa^*$ measurable is too high compared with the hope to
deal with $\bT\subseteq L_{\omega_1, \omega}$ (or any $L_{\kappa,
\omega}$) but seem quite small compared to the compact cardinal in
\cite{MaSh:285}. Model theoretically a compact cardinal ensure many cases
of amalgamation, whereas measurable cardinal ensure no maximal model. We
continue \cite{Sh:300}, \cite{MaSh:285}, \cite{KlSh:362}: try to
imitate \cite{MaSh:285}; 
a paralel line of research is \cite{Sh:394}. Earlier works are
\cite{Sh:48}, \cite{Sh:87a}, \cite{Sh:87b}; on the upward \L o\'s
conjecture, look at \cite{Sh:576} and \cite{Sh:600}.

On the situation with the upward direction and generlly more see \cite{Sh:576}.

This paper continues the tasks begun in \cite{KlSh:362}. We use the results
obtained therein to advance our knowledge of the categoricity spectrum of
theories in $L_{\kappa^*,\omega}$, when $\kappa^*$ is a measurable cardinal.

The main theorems are proved in section three; section one treats of types and
section two described some constructions.

The notation follows \cite{KlSh:362}, except in two important details:
we reserve
$\kappa^*$ for the fixed measurable cardinal and $\bT$ for the fixed
$\lambda$-categorical theory in $L_{\kappa^*,\omega}$ in a given vocabulary
$L$. $\kappa$ is any infinite cardinal and $T$ is usually some kind of
tree. To recap briefly: $\bT$ is a $\lambda$- categorical theory in
$L_{\kappa^*,\omega}$, $LS(\bT)\defeq
\kappa^*+|\bT|$, ${\cal K}=\langle K,\preceq_{\cal F}\rangle$
is the class of models of $\bT$, where ${\cal F}$ is a fragment
of $L_{\kappa^*,\omega}$ satisfying $\bT \subseteq{\cal F}$, $|{\cal F}|\leq
\kappa^*+|\bT|$, and for $M$, $N\in K$, $M\preceq_{\cal F} N$ means that
$M$ is an ${\cal F}$-elementary submodel of $N$.

The principal relevant results from \cite{KlSh:362} are: ${\cal K}_{<\lambda}$
has the amalgamation property (5.5 there) and every member of $K_{<\lambda}$
is nice (5.4 there).  But this assumption ($\bT$ categorical in
$\lambda$) or its consequences mentioned above will be mentioned
in theorems when used.

Let $(M_1,M_0)\preceq_{\cal F}(M_3,M_2)$ means $M_1\preceq_{\cal F} M_3$,
$M_0\preceq_{\cal F} M_2$.

$(I_1,I_2)$ is a Dededind cut of the linear order $I$
if  
$$
I=I_1\cup I_2, I_1\cap I_2=\emptyset, \forall x\in I_1\forall y\in I_2
(x<y),
$$
the two sided cofinality of $I$, $\dcf (I)$ is $(\cf I_1, \cf I_2^*)$
where $I^*_2$ is the order $I_2$ inverted.

Writing proofs we also consider their hopeful rule in the hopeful
classification 
theory. But we have been always carelul in stating the assumptions.

Note that \cite{KlSh:362} improve results of \cite{MaSh:285}; 
but they do not fully recapture the results on the compact case to the
measurable case, e.g. there the
results work for every $\lambda>\kappa^*$ whereas here we sometimes
need ``$\lambda$ above the Hanf number of omitting types'', say
$\beth_{(2^{LS(\bT)})^+}$. 

We thank Oren Kolman for writing and ordering notes from lectures on
the subject 
from spring 90 (you can see his style in the parts with good
language).

\section{Knowing the right types:}
The classical notion of type relates to the satisfaction of sets of formulas
in a model. We shall define a post-classical type (following
\cite{Sh:300}, \cite{Sh:h} which was followed by \cite{MaSh:285} but
niceness is involved) and use 
this to define notions of freeness and non-forking appropriate in the context
of a $\lambda$-categorical theory in $L_{\kappa^*, \omega}$.  The definitions
try to locate a notion which under the circumstances behave as in
\cite{Sh:c}. 

\begin{context}
\label{cont1.0}
{\rm
 $\bT\subseteq L_{\kappa^*, \omega}$ in the vocabulary $L$, $K=\{M:
 M\mbox{ a model of }\bT\}$, $\preceq_{\cF}$ as in the introduction.
$K_\mu=\{M\in K: \|M\|=\mu\}$,
$K_{<\kappa}=\bigcup\limits_{\mu<\kappa} K_\mu$, and ${\cal K}=(K,
\preceq_{\cF})$ and we stipulte $K_{<\kappa^*}=\emptyset$,
 e.g. $K_{<\kappa}=\bigcup\{K_\mu: \mu<\kappa\mbox{ but }\mu\geq
\kappa^*\}$. We let $LS(\cK)=|\cF|+\kappa^*$. 

Remember ``$M\in K$ is nice'' is defined in \cite{KlSh:362},
definitions 3.2, 1.8; nice implies being an amalgamation base in
$K_{<\lambda}$ (see \ref{conc3.3}).
}
\end{context}

\begin{definition}
\label{def1.1}
Suppose that $M\in K_{<\lambda}$ is a nice model of $\bT$.  Define a binary
relation, $E_M= E^{<\lambda}_M$, as follows:
\begin{description}
\item{\ } $(\bar a_1,N_1)E_M(\bar a_2, N_2)$ iff for $\ell=1,2,$ $N_{\ell}\in
K_{<\lambda}$ is nice and $M\preceq_{\cal F} N_{\ell}$, $\bar a_{\ell}\in
N_{\ell}$ (i.e. $\bar a_\ell$ a finite sequence of members of
$N_\ell$), and there exist a model $N$ and embeddings $h_{\ell}$ such that 
$M\preceq_{\cal F} N$, $h_{\ell}: N_{\ell} \mathop\to\limits_{\cal F} N$,
$id_M=h_1\restriction M=h_2\restriction M$, and $h_1(\bar a_1)=h_2(\bar a_2)$.
\end{description}
\end{definition}

\begin{fact}
\label{fact1.2}
\begin{enumerate}
\item $E_M$ is an equivalence relation.
\item Let $M\in K_{<\lambda}$, $M\preceq_{\cal F} N$, $\bar a\in N$,
and for $\ell=1, 2$, $N\cup \bar a \subseteq N_\ell\preceq_{\cal F}
M$, $\|N_\ell\| <\lambda$ then $(\bar a, N_1)E_M(\bar a, N_2)$
\end{enumerate}
\end{fact}

\Proof
1) To prove \ref{fact1.2}, let's look at transitivity.  

Suppose
$(\bar a_{\ell}, N_{\ell})E_M(\bar a_{\ell+1}, N_{\ell+1})$, $\ell=1,2$.
Thus there are models
$N^{\ell}$ and embeddings $h_0^{\ell}$, $h_1^{\ell}$ of $N_{\ell}$, $N_{\ell+1}$
over $M$ into $N^{\ell}$, with $h_0^{\ell}(\bar a_{\ell})=h_1^{\ell}(\bar
a_{\ell+1})$, $\ell=1,2$. W.l.o.g. $N^{\ell}\in K_{<\lambda}$ (by the Downward
Loewenheim Skolem Theorem). By assumption $N_2$ is nice,
hence by \cite[3.5]{KlSh:362}
is an amalgamation base for ${\cal K}_{<\lambda}$, i.e.  there is an amalgam
$N^*\in K_{<\lambda}$, and embeddings $g_{\ell} : N^{\ell}\stackrel{{\cal
F}}{\longrightarrow} N^*$, amalgamating $N^1$, $N^2$ over $N^2$ w.r.t
$h_1^1$, $h_0^2$. In other words, the following diagram commutes:
\medskip

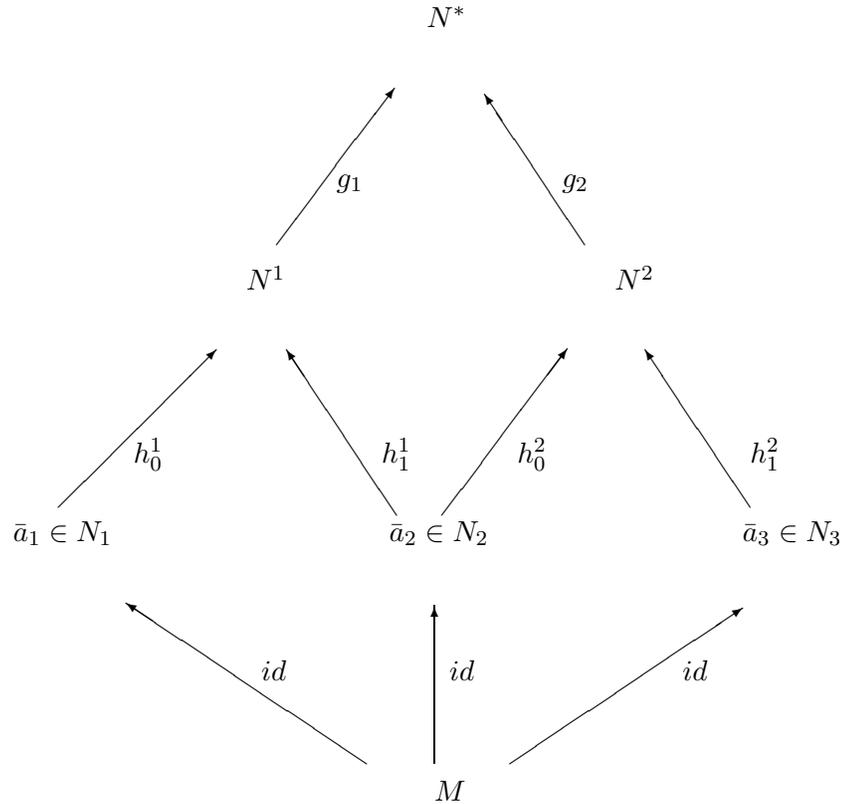
\unitlength=1.00mm
\special{em:linewidth 0.4pt}
\linethickness{0.4pt}
\begin{picture}(118.00,119.00)
\put(75.00,109.00){\makebox(0,0)[lb]{$N^*$}}
\put(51.00,74.00){\makebox(0,0)[lb]{$N^1$}}
\put(100.00,74.00){\makebox(0,0)[lb]{$N^2$}}
\put(55.00,80.00){\vector(3,4){15.67}}
\put(96.00,80.00){\vector(-2,3){13.33}}
\put(63.00,87.00){\makebox(0,0)[lb]{$g_1$}}
\put(93.00,87.00){\makebox(0,0)[lb]{$g_2$}}
\put(20.00,40.00){\makebox(0,0)[lb]{$\bar a_1\in N_1$}}
\put(26.00,45.00){\vector(1,1){21.00}}
\put(36.00,50.00){\makebox(0,0)[lb]{$h^1_0$}}
\put(70.00,40.00){\makebox(0,0)[lb]{$\bar a_2\in N_2$}}
\put(117.00,40.00){\makebox(0,0)[lb]{$\bar a_3\in N_3$}}
\put(71.00,44.00){\vector(-2,3){14.67}}
\put(69.00,50.00){\makebox(0,0)[lb]{$h^1_1$}}
\put(118.00,45.00){\vector(-2,3){14.00}}
\put(87.00,50.00){\makebox(0,0)[lb]{$h^2_0$}}
\put(118.00,50.00){\makebox(0,0)[lb]{$h^2_1$}}
\put(76.00,6.00){\makebox(0,0)[lb]{$M$}}
\put(76.00,11.00){\vector(0,1){21.00}}
\put(67.00,11.00){\vector(-3,2){32.00}}
\put(86.00,11.00){\vector(3,2){31.00}}
\put(53.00,22.00){\makebox(0,0)[lb]{$id$}}
\put(78.00,22.00){\makebox(0,0)[lb]{$id$}}
\put(109.00,22.00){\makebox(0,0)[lb]{$id$}}
\put(77.00,44.00){\vector(3,4){16.67}}
\end{picture}

\medskip

Just notice now that $N^*$, $g_1 h_0^1$, $g_2 h_1^2$ witness that
$(\bar a_1, N_1) E_M(\bar a_3, N_3)$, since:
$$g_1h_0^1(\bar a_1) = g_1(h_1^1(\bar a_2)) =
g_2h_o^2(\bar a_2) = g_2h_1^2(\bar a_3).
$$
\hfill\qed$_{\ref{fact1.2}}$

\begin{definition}
\label{def1.3}
Suppose that $M\in K_{<\lambda}$ is nice, $a\in N\in K_{<\lambda}$ and
$M\preceq_{\cal F} N$. 
Then
\begin{enumerate}
\item $\tp(a, M, N)$, the type of $a$ over $M$ in $N$,
is the $E_M$-equivalence
class of $(a, N)$,
$$
(a, N)/E_M=\{(b, N^1): (a, N) E_M(b, N^1)\}.
$$
We also say ``$a\in N$ realizes $p$''. If $\|N\|\geq \lambda$ define
$\tp(\bar a, M, N)$ by \ref{fact1.2}(2).

\item If $M'\preceq_{\cal F} M\in K_{<\lambda}$, $p\in S(M)$ (see
below) is $(a,N)/E_M$
then $p \restriction M'=(a,N)/ E_{M'}$.

\item If $LS(\bT)<\kappa\leq\mu\leq\lambda$, we call $M\in K_\mu$
$\kappa$-saturated
if for every nice $N\preceq_{\cal F} M$, $||N||<\kappa$ and $p\in S(N)$, some
$\bar a\in M$ realizes $p$ (in $M$) or at least for some $N^\prime$,
$N\preceq_{\cal F} N^\prime \preceq_{\cal F} M$, some $a^\prime\in
N^\prime$ realizes $p$ in $N^\prime$.

\item $S^m(N)=\{p: p=\tp(\bar a, N, N_1)$ for any $N_1$, $\bar a$
satisfying: $N\preceq_{\cal F} N_1$, $\| N_1\|\leq \|N\|+LS(\cK)$ and
$\bar a\in {}^m(N_1)\}$

$S^{<\omega}(N)=\bigcup\limits_{m<\omega} S^m(N)$.

\item $\bT$ is $\mu$-stable if $N\in K_{\leq \mu} \imply |S(N)|\leq
\mu$.

\item We say $N$ is $\mu$-universal over $M$ when: $M\preceq_{\cF} N$,
$N\in K_\mu$ and if $M\preceq_{\cF} N'\in K_{\leq \mu}$ {\em then} there
is a $\preceq_{\cF}$-embedding of $N'$ into $N$ over $M$.

\item We say $N$ is $(\mu, \kappa)$-saturated over $M$ if there is an
increasing continuous sequence $\langle M_i: i<\kappa\rangle$ such that:
$M_0=M$, $N=\bigcup\limits_{i<\kappa} M_i$, $M_i\in K_{\mu}$ and
$M_{i+1}$ is $\mu$-universal over $M_i$.

\item We say $\cK$ (or $\bT$) is stable in $\mu$ if for every $M\in
K_\mu$, $M$ is nice and $|S(M)|\leq \mu$.
\end{enumerate}
\end{definition}

\begin{definition}
\label{def1.4}
We shall write $\nonforkin{M_1}{M_2}_{M_0}^{M_3}$ to mean: $M_0\preceq_{\cF}
M_1\preceq_{\cF} M_3$, $M_0\preceq_{\cF} M_2 \preceq_{\cF} M_3$
and there exist suitable operation $(I, D, G)$ and an embedding $h:
M_3\stackrel{\cF}{\longrightarrow} \Op(M_1, I, D, G)$ such that
$h\restriction M_1 =\id_{M_1}$
and $\Rang(h\restriction M_2)\subseteq \Op(M_0, I, D, G)$ (remember that
$\Op(M, I, D, G)$ is the limit ultrapower of $M$ w.r.t. $(I, D, G)$; see
\cite[1.7.4]{KlSh:362}). We say that $M_1$, $M_2$ do not fork in
$M_3$ over $M_0$ if
$$\nonforkin{M_1}{M_2}_{M_0}^{M_3}.$$
If
$$\nonforkin{M_1}{M_2}_{M_0}^{M_3}$$
does not hold, we'll write
$$M_1\doesforkin_{M_0}^{M_3}M_2$$
and say that $M_1$, $M_2$ fork in $M_3$ over $M_0$.
\end{definition}

\begin{theorem}
\label{th1.5}

Suppose that $\nonforkin{M_1}{M_2}_{M_0}^{M_3}$ and
$M_2\doesforkin_{M_0}^{M_3}M_1$ (failure of $\nonfork{}{}_{}$-symmetry)
and $M_0\prenice M_3$.

Let $\mu=\kappa^*+|{\bT}|+||M_2||+||M_1||$. Then for every linear order
$(I, <)$ there exists an Ehrenfeucht-Mostowski model $N=EM(I, \Phi)$ with
$\mu$ (individual) constants $\{\tau_i^0 : i<\mu\}$ and unary function symbols
$\{\tau_i^1(x_i) : i<\mu\}$, $\{\tau_i^2(x_i) : i<\mu\}$ such that, for
$M=(N\restriction L) \restriction \{\tau_i^0 : i<\mu\}$ (i.e. $M$ is a submodel of $N$
with the same vocabulary as $\bT$ and universe $\{\tau^0_i:i<\mu\}$ i.e. the
set of interpretations of these individual constants and for every $t\in I$,
$\ell=1,2,$
$$
M_t^\ell=(N\restriction L) \restriction \{
\tau_i^\ell(x_i):i<\mu\},
$$
one has $M\preceq_{\cF} N$,
$M_t^\ell\preceq_{\cF} N$ and for $s\not=t\in I$, $t<s$ iff
$\nonforkin{M_t^1}{M_s^2}_M^N.$
\end{theorem}

\noindent {\bf Remark:} Note $M_0\prenice M_3$ is automatic in the interesting
case since $M_0\in K_{<\lambda}$ and every element of $K_{< \lambda}$ is nice
by \cite[5.4]{KlSh:362}.

\noindent On the opertions see \cite{KlSh:362}.
\medskip

\Proof W.l.o.g. $\|M_3\|=\mu$. Add Skolem functions to $M_3$. We know that $M_0
\prenice M_3$. So there is $\Op^1$ such that $M_0\preceq_{\cF}
M_1\preceq_{\cal F} \Op^1 (M_0)$ and $\Op^2$ such that $M_1\preceq_{\cF}
M_3\preceq_{\cF} \Op^2(M_1)$, $M_2\preceq_{\cF} \Op^2(M_0)$. Let
$\Op=\Op^2\circ \Op^1$.
For each $t\in I$, let $\Op_t=\Op$. Let $N$ be
the iterated ultrapower of $M_0$ w.r.t.  $\langle Op_t : t\in I\rangle$. For
each $t\in I$, there is a canonical ${\cF}$-elementary embedding $F_t :
Op_t(M_0)\stackrel{{\cF}}{\longrightarrow} N$.
Let
$M=M_0$, and $M_t^\ell = F_t(M_\ell)$ for $\ell=1,2, t\in I$.

For each $t<s$, we can let $M^+_s=\langle \Op_v: v<s\rangle (M_0)$, so
$M_0\preceq_{\cF} M^+_t\preceq_{\cF} M^+_s\preceq_{\cF} \Op^1(M^+_t)$ and
we can extend $F_t\restriction M_1$ to an embedding of $\Op^2(M_1)$ into
$\Op^2_t(\Op^1_t(M^+_t))$, so $(F_t\restriction M_1) \cup (F_s\restriction
M_2)$ can be extended to a $\preceq_{\cF}$-embedding of $M_3$ into
$N$.
From the definition of the iterated ultrapower it follows that for
$s\neq t\in I$, $t<s$ implies  $\nonforkin{M_t^1}{M_s^2}_{M_0}^{N}$
and on the other hand 
by the assumption it follows that 
if $s$, $t\in I$, $s<t$ then $M^1_t
\doesforkin_{M_0}^{N} M^2_s$.
\hfill\qed$_{\ref{th1.5}}$

\begin{corollary}
\label{corol1.6}
Assume $\bT$ categorical in $\lambda$ or just $I(\lambda,
\bT)<2^\lambda$. Then 
$\bigcup\limits_{\mu^+<\lambda}K_{\mu}$ obeys $\nonfork{}{}_{}$-symmetry,
i.e.:  if $\nonforkin{M_1}{M_2}_{M_0}^{M_3}$ holds for $M_0$, $M_1$,
$M_2$, $M_3 \in
\bigcup\limits_{\mu^+<\lambda} K_{\mu}$, then
$\nonforkin{M_2}{M_1}_{M_0}^{M_3}$ holds.
\end{corollary}

\Proof If $\mu^+<\lambda$, $\nonforkin{M_1}{M_2}_{M_0}^{M_3}$ and
$M_2\doesforkin_{M_0}^{M_3}M_2$, then 
theorem \ref{th1.5} gives the assumptions of 
the results at the end of section
three in \cite{Sh:300}, III (or better \cite{Sh:e}, III, \S 3), yield
a contradiction to the $\lambda$-categoricity of $\bT$ and even
$2^\lambda$ pairwise non isomorphic models.
\hfill\qed$_{\ref{corol1.6}}$

It may be helpful, though somewhat vague, to add the remark that
$\nonfork{}{}_{}$-asymmetry enables one to define order and to build many
complicated models; so \ref{corol1.6} removes a potential
obstacle to a categoricity
theorem.

\begin{definition}
\label{def1.7}
Let $A$ be a set. We write $\nonforkin{M_1}{A}_{M_0}^{M_3}$ (where
$A\subseteq M_3$, $M_0\preceq_{\cF} M_1\preceq_{\cF} M_3$)
to mean that there exist $M_2$, $M'_3$ such that $A\subseteq |M_2|$,
 $M_3\preceq_{\cF} M'_3$
and $\nonforkin{M_1}{M_2}_{M_0}^{M'_3}$.  In this situation we say
that $A/M_1=\tp(A, M_1, M_3)$
does not fork over $M_0$ in $M_3$.

We'll write $\nonforkin{M_1}{a}_{M_0}^{M_3}$ to mean
$\nonforkin{M_1}{\{ a \}}_{M_0}^{M_3}$,
we then say $\tp(a, M_1, M_3)$ do not
fork on $M_0$.

We write $\nonforkin {A_1}{A_2}_{M_0}^{M_3}$ if for some $M_3$,
$M_3\preceq_{\cF} M_3'\in K_{<\lambda}$, and for some $M'_1$,
$A_2\subseteq M'_1\preceq_{\cF} M'_3$, and
$\nonforkin{M'_1}{A_2}_{M_0}^{M'_3}$.

\end{definition}

\begin{remark}
\label{remark1.8}
{\rm 
\begin{enumerate}
\item  Of particular importance is the case where $A$ is finite. Let us
explain the reason. We wish to prove a result of the form:
\begin{description}
\item[$(*)$] if $\langle M_i : i\le \delta+1\rangle$ is a continuous
$\prec_{\cF}$-chain and $a\in M_{\delta}$, then there is $i<\delta$ such that
$\nonforkin{M_\delta}{a}_{M_i}^{M_{\delta+1}}$.
\end{description}
This says roughly that the type $\tp(a, M_\delta, M_{\delta+1})$
is definable over a
finite set (or at least in some sense has finite characters). In general the
former
relation is not obtained. However its properties are correct.  Hence it will
be possible to define the rank of a over $M_0$, $\rrk(a, M_0)$, as an ordinal,
so that for large enough $M_3$, if
$M_1 \doesforkin_{M_0}^{M_3} a$, then
$\rrk(a, M_1)<\rrk(a, M_0)$.

\item If $A$ is an infinite set, then we cannot prove $(*)$, in general. For
example, suppose that $\langle M_i : i\le \omega\rangle$ is (strictly)
increasing continuous, $a_i\in (M_{i+1}\setminus M_i)$ and $A=\{a_i :
i<\omega\}$. Then for every $i< \omega$,
$\big(\bigcup\limits_{j<\omega}M_j\big)\doesforkin_{M_i}^{M_{\omega}}A$.  Still we can
restrict ourselves to $\delta$ of cofinality $>|A|$.

\item Notice that quite generally speaking,
      $\nonforkin{N_1}{N_2}_{N_0}^{N_1}$ implies that $N_1\cap N_2=N_0$.
\end{enumerate}
}
\end{remark}

\begin{definition}
\label{def1.9}
We define
\begin{description}
\item{\ }
$\boldk_\mu({\bT})=\boldk_\mu({\cK})=\{\kappa:
\cf(\kappa)=\kappa\leq\mu$ and there exist a continuous $\prec_{\cF}$-chain
$\langle M_i: i\leq\kappa+1\rangle\subseteq K_{\leq\mu}$ and $a\in
M_{\kappa+1}$ such that for all $i<\kappa$, $a/M_{\kappa}$ forks over $M_i$ in
$M_{\kappa+1}\}$.
\end{description}
I.e. for $\kappa\in \boldk_\mu({\bT})$ there are
$\langle M_i\in K_{\leq\mu} : i\leq\kappa+1\rangle$
and $a \in M_{\kappa+1}$ such
that $i<\kappa \imply
M_{\kappa}\doesforkin_{M_i}^{M_{\kappa+1}}a$.
\end{definition}

\begin{example}
\label{exam1.10}
{\rm Fix $\mu$ and $\alpha\leq\mu$.
Let $({\ }^\mu\omega, E_{\beta})_{\beta<\alpha}$ be the structure
with universe
$$
{ }^\mu\omega=\{\eta : \eta \mbox{ is a function from }\mu\mbox{ to }\omega\},
$$
$\eta E_{\beta} \nu$ iff $\eta \restriction \beta =\nu \restriction \beta$.
Let
$\bT=Th({}^\mu\omega, E_{\beta})_{\beta<\alpha}$.  Then
$\boldk_{\mu}(\bT)=\{\kappa : \cf(\kappa) =\kappa \leq\alpha \}$.}
\end{example}
Why? If $\cf(\kappa)=\kappa \leq\alpha$, then there are $M_i (i\leq\kappa+1)$,
$a\in M_{\kappa+1}$ and $a_i\in (M_{i+1}\setminus M_i)$ such that $a_i/
E_{i+1}\notin M_i$ (that's to say, no element of $M_i$ is $E_{i+1}$-
equivalent to $a_i$) and $aE_ia_i$.

\begin{definition}
\label{def1.11}
The class ${\cK}=\langle K, \preceq_{\cF}\rangle$ is $\chi$-based iff
for every pair of continuous $\prec_{\cF}$-chains $\langle N_i\in
K_{\leq\chi} : i<\chi^+\rangle $, $\langle M_i\in K_{\leq\chi} : i <
\chi^+\rangle$, with $M_i\preceq_{\cF} N_i$, there is a club $C$ of
$\chi^+$ such that
$$(\forall i\in C)\left(\nonforkin{M_{i+1}}{N_i}_{M_i}^{N_{i+1}}\right).$$

Replacing $\chi^+$ by regular $\chi$ we write $(<\chi)$-based.  We say
synonymously that $\bT$ is $\chi$-based.
\end{definition}

\begin{definition}
\label{def1.12}
The class ${\cK}=\langle K, \preceq_{\cF}\rangle$ has continuous
non-forking in $(\mu, \kappa)$ iff
\begin{description}
\item[$(\alpha)$] {\em whenever} $\langle M_i\in K_{\leq\mu} : i\leq
\delta\rangle$ is an continuous $\prec_{\cF}$-chain, $|\delta|\leq\mu$,
$\cf(\delta)=\kappa$,
$$M_0\preceq_{\cF} N_0\preceq_{\cF} N^*,\ M_{\delta}\preceq_{\cF}
N^*\mbox{ and }$$
$$(\forall i<\delta)\left(\nonforkin{M_i}{N_0}_{M_0}^{N^*}\right),$$
{\em then}
$\nonforkin{M_{\delta}}{N_0}_{M_0}^{N^*}$;

\item[$(\beta)$] {\em whenever} $\langle M_i\in K_{\leq\mu} :
i\leq\delta+1\rangle$, $\langle N_i\in K_{\leq\mu} : i\leq\delta+1\rangle$ are
continuous $\prec_{\cF}$-chains, $M_i\preceq_{\cF} N_i$,
$|\delta|\leq\mu$, $\cf(\delta)=\kappa$ and
$$(\forall
i<\delta)\left(\nonforkin{M_{\delta+1}}{N_i}_{M_i}^{N_{\delta+1}}\right),$$
{\em then} $\nonforkin{M_{\delta+1}}{N_{\delta}}_{M_{\delta}}^{N_{\delta+1}}$.
\end{description}
Again we'll mean the same thing by saying that $\bT$ has continuous
non-forking in $(\mu,\kappa)$.
\end{definition}

Our next goal is to show that if $\bT$ fails to possess these features for
some $\mu<\lambda$ such that $\mu\geq \kappa + LS(\cK)$, then $\bT$ has many models in $\lambda$.


Let us quote in this context a further important result from
\cite{Sh:300}, II, 2.12:

\begin{theorem}
\label{th1.13}
Assume $\bT$ be a $\lambda$-categorical theory, or just $K_{<\lambda}$
has amalgamation and every $N\in K_{<\lambda}$ is nice.
\begin{enumerate}
\item Let $LS(\bT)<\mu\leq\lambda$, $M\in K_\mu$. Then TFAE:
\begin{description}
\item[(A)] $M$ is universal-homogeneous: if $N\preceq_{\cF}M$, $\|N\|<\mu$,
$N\preceq_{\cF}N'\in K_{<\mu}$, {\em then} there is an ${\cF}$-elementary
embedding $g: N'\stackrel{{\cF}}{\longrightarrow}M$ such that
$g\restriction N=\id_N$.
\item[(B)] if $N\preceq_{\cF} M$, $||N||<\mu$ and $p\in S(N)$,
{\em then} $p$ is realized in $M$ i.e. $N$ is saturated.
\end{description}
\item $M$ as in (A) or (B) is unique for fixed $\bT$, $\mu$.
\item
Any two $(\mu, \kappa)$-saturated models are isomorphic (see
\ref{def1.3}(7)). 
\end{enumerate}
\end{theorem}

\Proof 1), 2) See \cite{Sh:300}, II 3.10.

3) Easy.
\hfill\qed$_{\ref{th1.13}}$

\begin{claim}
\label{claim1.14}
Assume $\bT$ is $\lambda$-categorical or just $\cK_{<\lambda}$ has
amalgamation. 
\begin{enumerate}
\item If $LS(\bT)\leq\mu<\lambda$, $N_0\preceq_{\cF} N_1$ are in $K_\mu$ then
TFAE
\begin{description}
\item[(A)] $N_1$ is $(\mu,\mu)$-saturated over $N_0$

\item[(B)] there is a $\preceq_{\cF}$-increasing continuous $\langle
M_i:i\leq\mu\times \mu\rangle$, such that:

$M_\mu=N_1$, $M_0=N$ and every $p\in S(M_i)$ is realized in $M_{i+1}$
\end{description}
\item Also TFAE for $\kappa=\cf (\kappa)\leq \mu^+$
\begin{description}
\item[(A)$_\kappa$] $N_1$ is $(\mu, \kappa)$-saturated over $N_0$

\item[(B)$_\kappa$] there is a $\preceq_{\cF}$-increasing continuous
$\langle M_i: i\leq \mu\times \kappa\rangle$ with $M_{\mu\times
\kappa}=N_1$, $M_0=N$ and every $p\in S(M_i)$ is realized in $M_{i+1}$
\end{description}
\item If $\cK$ is stable in $\mu$, $\mu\geq LS(\cK)$,
$\kappa=\cf(\kappa)\leq \mu^+$ {\em then} there is a $(\mu, \kappa)$-saturated
model.
\end{enumerate}
\end{claim}

\Proof
1) See \cite{Sh:300}, II 3.10

2) Same proof.

3) Straight.
\hfill\qed$_{\ref{claim1.14}}$

\begin{claim}
\label{claim1.15}
($\bT$ categorical in $\lambda$)
\begin{enumerate}
\item Any $M\in K_\lambda$ is saturated.
\item Every $N\in K_{<\lambda}$ is nice.
\item $K_{<\lambda}$ has $\preceq_{\cF}$-amalgamation.
\item $\cK$ is stable in $\mu$ for $\mu\in [LS(\bT), \lambda)$.
\end{enumerate}
\end{claim}

\Proof
1) By the proof of \cite[5.4]{KlSh:362} (for $\lambda$-regular easier).

2) \cite[5.4]{KlSh:362}.

3) \cite[5.5]{KlSh:362}.

4) See \cite{KlSh:362}. 
\hfill$\qed_{\ref{claim1.15}}$

\begin{intermediatecorollary}
\label{inter1.15A}
\begin{enumerate}
\item Suppose that {\bf T} is $\lambda$-categorical. If $\mu<\lambda$,
$\mu>LS(\hbox{\bf T})$ and {\bf T} is not $\mu$-categorical, then there is an
unsaturated model $M \in K_\mu$.

\item It now follows that if we show that the existence of an
unsaturated model in $K_\mu$ implies that of an unsaturated model in
$K_\lambda$, then $\lambda$-categoricity of {\bf T} implies $\mu$-categoricity
of {\bf T}. 

\end{enumerate}
\end{intermediatecorollary}

\begin{conclusion}
\label{conc1.16}
[$\bT$ categorical in $\lambda$]
If $I$ is a linear order, $I=I_1+I_2$, $|I|<\lambda$ and $J=I_1+\omega
+I_2$ {\em then} every $p\in S(EM(I))$ is realized in $EM(J)$.
\end{conclusion}

\Proof $EM(I_1+\lambda+I_2)$ is in $K_\lambda$ hence is saturated hence every
$p\in S(EM(I))$ is realized in it, say by $a_p$, for some finite
$w_p\subseteq \lambda$ we have $a_p\in EM(J_1+w_p+I_2))$, now we
indiscernibility.
\hfill\qed$_{\ref{conc1.16}}$

\begin{remark}
\label{rem1.16A} 
{\rm By changing $\Phi$ we can replace ``$\omega$''
by ``$1$''.}
\end{remark}

\begin{conclusion}
\label{conc1.17}
[$\bT$ categorical in $\lambda$]
\begin{enumerate}
\item If $J=\bigcup\limits_{\alpha<\mu}I_\alpha$,
$|J|=\mu\in[LS(\bT),\lambda]$, 
$I_\alpha$ increasing continuous, for each $\alpha$ some Dedekind cut of
$I_\alpha$ is realized by infinitely many members of $I_{\alpha+1}\setminus
I_\alpha$ {\em then} $EM(J)$ is $(\mu,\mu)$-saturated over $EM(I_0)$.
\item If $\Phi$ is ``corrected'' as in \ref{rem1.16A}, $I_0\subseteq
J$, $|J\setminus I_0|= |J|=\mu$, $\mu\in [LS(\bT), \lambda]$, {\em
then} $EM(J)$ is $(\mu, \mu)$-saturated over $EM(I_0)$ moreover for
any $\kappa=\cf(\kappa)\leq \mu$ it is $(\mu, \kappa)$-saturated.
\end{enumerate}
\end{conclusion}

\Proof By \ref{conc1.16}+\ref{claim1.14}.

\begin{claim}
\label{claim1.18}
\begin{enumerate}
\item Suppose $\langle N_i^\ell:i\leq\alpha\rangle$ is
$\prenice$-increasing continuous for $\ell=1,2$, $N_i^1\preceq_{\cF}
N_i^2\in K_{<\lambda}$ and
$\nonforkin{N_i^2}{N_{i+1}^1}_{N_i^1}^{N_{i+1}^2}$ for each $i<\alpha$, then
$\nonforkin{N_0^2}{N_\alpha^1}_{N_0^1}^{N_\alpha^2}$.

\item The monotonicity properties of $\nonforkin{}{}_{}^{}$, i.e.: if
$\nonforkin{M_1}{M_2}_{M_0}^{M_3}$ and for some operation $\Op$ and
moduls $M'_1$, $M'_2$, $M'_3$ we have
$M_3\preceq_{\cF} M'_3\preceq \Op(M_3)$ and $M_0 \preceq_{\cF} M'_1
\preceq_{\cF} M_1$ and $M_0\preceq_{\cF} M'_2\preceq_{\cF} M_2$ then
$\nonforkin{M'_1}{M'_2}_{M_0}^{M'_3}$.

\item If $\nonforkin{M_1}{A}_{M_0}^{M_3}$ and $M_0\preceq_{\cF}
M'_0\preceq_{\cF} M'_1\preceq_{\cF} M'_3 \preceq_{\cF} M''_3$ and
$M_3\preceq_{\cF} M''_3$ and $A'\subseteq A$ then
$\nonforkin{M'_1}{A'}_{M'_1}^{M'_3}$. 
\end{enumerate}
\end{claim}

\Proof Use \cite[1.11]{KlSh:362}. 

\begin{claim}
\label{claim1.19}
[$\bT$ is $\lambda$-categorical]
If $M_0\prenice M_1, M_2$ are in $K_{<\lambda}$ {\em then} we can find $M_4\in
K_{<\lambda}$, $M_0\preceq_{\cF} M_4$ and $\preceq_{\cF}$-embeddings
$f_1$, $f_2$ of $M_1$, $M_2$ respectively into $M_4$ such that
\begin{description}
\item[($\alpha$)] $\nonforkin{f_1(M_1)}{f_2(M_2)}_{M_0}^{M_4}$ and
\item[($\beta$)] $\nonforkin{f_2(M_2)}{f_1(M_1)}_{M_0}^{M_4}$.
\end{description}
\end{claim}

\begin{remark}
\label{rem1.19A}
{\rm Note \ref{corol1.6} deal only with models in $\bigcup\{K_\mu:
\mu^+<\lambda\}$ hence $(\beta)$ is not totally redundant.} 
\end{remark}

\Proof If we want to get ($\alpha$) only, use operation $\Op$ such that
$\Op(M_0)$ has cardinality $\geq\lambda$, choose $N\preceq_{\cF} \Op(M_0)$,
$\|N\|=\lambda$, hence $N$ is saturated hence we can find
a $\preceq_{\cF}$-embedding $f_2: M_2\to N$, let
$N_1=\Op(M_1)$, so $N\preceq_{\cF} \Op(M_0)\preceq_{\cF}\Op(M_1)=N_1$,
and choose $M_4\prec N_1$,
$M_4\in K_\mu$ such that $M_1\cup\Rang f_2\subseteq N$.

By ``every $N\in K_\lambda$ is saturated'' there are an operation
$\Op$ and $N\in 
K_\lambda$ such that  $M_0\preceq_{\cF} N\preceq_{\cF} \Op (M_0)$
hence there are 
$M_0^+$, $M_1^+$, $M_2^+$ in $K_{<\lambda}$ such that:
\begin{description}
\item[$(*)_0$] $(M_1^+,M_0^+)\preceq_{\cF} \Op(M_1,M_0)$,
$(M^+_2,M_0^+)\preceq_{\cF} \Op(M_2,M_0)$ and $M_0^+$ has the form
$EM(I_0)$, $I_0$ a linear order with $|I_0|$ Dedekind cuts with cofinality
$(\kappa^*,\kappa^*)$. [Note that by \ref{conc1.17}(2) if $LS(\bT)\leq
|I_0|\leq\lambda$ 
then $EM(I_0)$ is saturated and $N$ is saturated, clearly there is
$I_0$ as required.]
\end{description}
Hence we can find $I_1,I_2,I_3$ such that: $I_o\stackrel{\rm def}{=}
 I\subseteq I_1\subseteq I_3$ $I_0\subseteq I_2\subseteq I_3$,
$I_1\cap I_2=I$, no $t_1\in
I_1\setminus I_0$, $t_2\in I_2\setminus I_0$ realize the same Dedekind cut of
$I$, and every $t\in T_3\setminus I_0$ realize a cut of $I$ with cofinality
$(\kappa^*,\kappa^*)$. Hence $I_0\mathop\subseteq_{\nice} I_\ell$
$(\ell\geq3)$, moreover $\nonforkin{I_1}{I_2}_{I_0}^{I_3}$ and
$\nonforkin{I_2}{I_1}_{I_0}^{I_3}$ hence
\begin{description}
\item[$(*)_1$] $\nonforkin{EM(I_1)}{EM(I_2)}_{EM(I_0)}^{EM(I_3)}$ and
$\nonforkin{EM(I_2)}{EM(I_1)}_{EM(I_0)}^{EM(I_3)}$.
\end{description}
Also by \ref{conc1.17}(2), wlog ($\ell=1,2$) $M^+_\ell\preceq_{\cF}
EM(I_\ell)$. 
\par \noindent
So
\begin{description}
\item[$(*)_2$] $\nonforkin{M^+_1}{M_2^+}_{M_0^+}^{EM(I_3)}$ and
$\nonforkin{M_1^+}{M_1^+}_{M_0^+}^{EM(I_3)}$.
\end{description}

By $(*)_1+(*)_2$ and \ref{claim1.18}(1) (for $\alpha=2$) we get the
conclusions. 
\hfill\qed$_{\ref{claim1.19}}$

\begin{claim}
\label{claim1.20}
[$\bT$ is  $\lambda$-categorical]

\begin{enumerate}
\item If $\nonforkin{M_1^\ell}{M_2^\ell}_{M_0^\ell}^{M_3^\ell}$
for $\ell=1,2$, $M_\ell^3\in K_{<\lambda}$ moreover $\|M_3^\ell\|^+<\lambda$
and $f_k$ an isomorphism from $M^1_k$ onto $M^2_k$ for $k=0,1,2$
such that $f_0\subseteq
f_1$, $f_0\subseteq f_2$ {\em then} there is $M$, $M_2^3\preceq_{\cF} M\in
K_{<\lambda}$, $\|M\|=\|M_3^1\|+\|M_3^2\|$ and a
$\preceq_{\cal F}$-embedding $f$
of $M_3^1$ into $M_3^2$ extending $f_1$ and $f_2$.

\item
Assume $\nonforkin{M_1^\ell}{M_2^\ell}_{M_0^\ell}^{M_3^\ell}$ for
$\ell=1,2$ and $A_2^\ell\subseteq M_2^\ell\preceq M_3^\ell$, and
$M_3^\ell\in K_{<\lambda}$ moreover $\|M_3^\ell\|^+<\lambda$, and $f_k$
is an isomorphism from $M_k^1$ onto $M^2_k$ for $k=0,1,2$ such that
$f_0\subseteq f_1$ and $f_0\subseteq f_2$ and $f_2$ maps $A^1_2$ onto
$A^2_2$ {\em then} there is $M$, $M^3_2\preceq_{\cF} M\in K_{<\lambda}$
such that $\|M\|=\|M^1_3\|+\|M^2_3\|$ and a $\preceq_{\cF}$-embedding
$f$ of $M^1_3$ into $M^2_3$ extending $f_1$ and $f_2\restriction A^1_2$.

\item If for $\ell=1,2$ we have $p_\ell\in S(N)$ does not fork over $M$ (see
Definition \ref{def1.7}),
$M\preceq_{\cF} N\in K_{\mu}$, $\mu^+<\lambda$ and $p_1\restriction M=
p_2\restriction M$ {\em then} $p_1=p_2$
\end{enumerate}
\end{claim}

\begin{remark}
\label{rem1.20A}
{\rm
\begin{enumerate}
\item This is uniqueness of non forking amalgamation.

\item The requirement is $\|M^\ell_3\|^+<\lambda$ rather than
$\|M^\ell_3\|<\lambda$ only because of the use of symmetry, i.e.
\ref{corol1.6}.
\end{enumerate}
}
\end{remark}

\Proof Wlog $f_0=\id$, $M_0^1=M_0^2$ call it $M_0$
and $f_1=\id_{M^1_1}$, $M_1^1=M_1^2$ call it
$M_1$. For some operation $\Op_{\ell}$ we have  $(M_3^\ell,
M_2^\ell)\preceq_{\cF} 
\Op_\ell(M_1^\ell,M_0^\ell)$. Let $\Op=\Op_1\circ \Op_2$, so
$M_3^\ell\preceq_{\cF} \Op(M_1)$, $M_2^\ell\preceq_{\cF} \Op(M_0)$.
W.l.o.g. $\|\Op(M_0)\|\geq\lambda$ and $\|\Op(M_1)\|\geq\lambda$,
so there is $N_0$,
$\bigcup\limits_{\ell=1}^2 M_2^\ell\subseteq N_0\preceq_{\cF} \Op(M_0)$,
such that
$\|N_0\|=\lambda$, hence $N_0$ is saturated hence there is an automorphism
$g_0$ of $N_0$ such that $g_0\restriction M_2^1=f_2$ (so $g_0\restriction
M_0=\id_{M_0})$. So there is $N_2$, $\bigcup\limits_{\ell=1}^2M_2^\ell\subseteq
N_2\preceq_{\cF} N_0$, $\|N_2\|<\lambda$, $N_2$ closed under $g_0$,
$g_0^{-1}$. Now there is $N_3$, $N_0\cup M_1\subseteq N_3\preceq_{\cF}
\Op(M_1)$, $N_3\in K_\lambda$, hence $N_3$ is saturated. 
So $\nonforkin{M_1}{N_2}_{M_0}^{N_3}$ hence
$\nonforkin{N_2}{M_1}_{M_0}^{N_3}$ (by symmetry i.e. \ref{corol1.6})
hence for some $N'_3$, $N_3\preceq_{\cF} 
N'_3\in K_{<\lambda}$, some automorphism $g_1$ of $N_3$ extend
$(g_0\restriction N_1)\cup \id_{M_1}$.  [Why? for some $\Op'$, $(N_3,
M_1)\preceq_{\cF} \Op'(N_1, M_0)$ 
and $\Op'(N_1)$, $\Op'(g_0\restriction N_1)$ are as required except having too
large cardinality, but this can be rectified.]

Clearly we are done.

2) Similarly.

3) Follows.
\hfill\qed$_{\ref{claim1.20}}$

\section{Various constructions}

In this section we'll attempt to describe some constructions of models of
$\bT$
relating to the situations in \ref{def1.11} and \ref{def1.12}, i.e. we
want to prove there are ``many complicated'' models of $\bT$ when
$\bT$ is ``on unstable side'' of Def.\ref{def1.11} or
Def.\ref{def1.12}.  May we 
suggest that on a 
first reading the reader be content with the perusal of \ref{const2.1}
and \ref{const2.2}, leaving 
the heavier work of \ref{improv2.2.4} until after section three which
contains the 
model-theoretic fruits of the paper.  The construction should be meaningful
for the classification problem.

What we actually need are \ref{moddata2.2.1}, \ref{claim2.2.2},
\ref{fact2.2.3} 
\begin{construction}{First try}
\label{const2.1}
\  
\end{construction}

\begin{data}
\label{data2.2.1}
{\em
Suppose that $\langle M_i\in K_{\leq\mu} : i\leq \kappa+1\rangle$ is an
continuous $\prec_{\cF}$-chain of models of $\bT$, $\mu<\lambda$;
$T$ is a non empty subset of $( ^{\kappa+1\geq}Ord)$ and
\begin{description}
\item[(i)] $T$ is closed under initial segments, i.e. if $\eta\in T$ and
$\nu\triangleleft \eta$, then $\nu\in T$,
\item[(ii)] if $\eta \in T$ and $\lg(\eta)=\kappa$ then
$\eta^\wedge\langle 0\rangle\in T$ and for
all $i$, $\eta^\wedge \langle 1+i\rangle\not\in T$.
\end{description}
Let $\lim_\kappa(T)=\{\eta :\lg(\eta)=\kappa\mbox{ and }
\mathop\wedge\limits_{i<\kappa}
(\eta\restriction i\in T)\}$.  Let $\{\eta_i : i<i^*\}$ be an enumeration of
$T$ such that if $\eta_i\triangleleft \eta_j$ ($\eta_i$ is an initial segment
of $\eta_j$), then $i<j$, and if $\eta_i=\nu^\land\langle \alpha \rangle$,
$\eta_j=\nu^\land\langle\beta\rangle$, $\alpha<\beta$, then $i<j$. For
simplicity $i^*$ is a limit ordinal.
}
\end{data}

\begin{ftry}
\label{details2.2.2}
{\em
From the data of \ref{data2.2.1} we shall build a model $N^*$
with Skolem functions,
$N^*\restriction L\in K$, and for $\eta\in T$, $M^*_\eta\subseteq N^*$,
$f_\eta: M_{lg(\eta)}
\mathop{\mapsto}\limits_{\cF}^{\mbox{onto}}M^*_\eta\restriction L$
such that if $\eta_i\triangleleft \eta_j$, then
$f_{\eta_i}\subseteq f_{\eta_j}$,
and $M^*_{\eta_i}\preceq_{{\cF}^{sk}}M^*_{\eta_j}$,
where ${\cF}^{sk}\supseteq$ $\bT^{sk}$ is a fragment of
$(L^{sk})_{\kappa^*,\omega}$.

Let $M^*_i=Sk(M_i)$ be a Skolemization of $M_i$ for $\cF$, increasing
($\subseteq$) with $i$ i.e. for every formula $(\exists y)\varphi(y,
\bar x)\in {\cal F}$ we choose a function $F^{M_i}_{\varphi(y, \bar
x)}$ from $M_i$ to $M_i$, with $\lg(\bar x)$-places such that 
$$
M_i \models (\exists y)\varphi(y, \bar a) \rightarrow
\varphi(F^{M_i}_{\varphi(y, \bar x)}(\bar a), \bar a)\mbox{ such that
}j<i \Rightarrow F^{M_i}_{\varphi(y, \bar y)}\restriction M_j=
F^{M_j}_{\varphi(y, \bar x)}.
$$
Note: we do not require even $M^*_i\prec M^*_{i+1}$.
}
\end{ftry}

To achieve this, let us define by induction on $i\leq i^*$, $N_i^*$,
$M_{\eta_i}^*$ and $f_{\eta_i}$. W.l.o.g. $\eta_0=\langle\ \rangle$. Let
$N_0^*=M^*_{\eta_0}= Sk(M_0)$, the Skolemization of $M_0$, $f_{\langle
\ \rangle}=\id_{M_0}$. If $i$ is a limit ordinal, let
$N_i^*=\bigcup\limits_{j<i}N_j^*$.  If $i$ is a successor ordinal and
$\lg(\eta_i)=\alpha+1$, then letting $\eta_j=\eta_i\restriction\alpha$, note
that $\eta_j\triangleleft \eta_i$ so $j<i$ and so $M^*_{\eta_j}$ and
$f_{\eta_j}$ are defined. We are assuming
$M_\alpha\prenice M_{\alpha+1}$ hence,
there is an operator $\Op=\Op_\alpha$ such that $M_\alpha\prenice
\Op(M_\alpha)$. Let 
$N^*_i=\Op(N^*_{i-1})$, let
$\Op(N^*_{i-1},M_\alpha,f_{\eta_j})=(N^*_{i-1},\Op(M_\alpha),(\Op(f_{n_i}))$,
and
let $f_{\eta_i}=\Op(f_{\eta_j})\restriction M_{\lg(\eta_i)}$ and
$M^*_{\eta_i}=\Rang(f_{\eta_i})$. (We can replace $N^*_{i+1}$ by any
$N'$ such that $N^*_i\cup M^*_{\eta_i}\subseteq N'\preceq_{\cF} N^*_{i+1}$
so preserving
$|N^*_i|\leq\mu+|i|$.  Finally, let $N^*=\bigcup\limits_{i<i^*}N_i^*$.  We are
left with the case $i$ successor ordinal, $\lg (\eta_i)$ a limit ordinal; we let
$N^*_i=N^*_{i+1}$, $M^*_{\eta_i}=\bigcup\limits_{\nu\triangleleft
\eta_i}M^*_\nu$ and $f_{\eta_j}=\bigcup\limits_{\nu\triangleleft \eta_j}f_\nu$.

\noindent{\bf Explanation:} In order to use this construction to prove
non-structure results, we intend to use property: for every
$\eta\in\lim_\kappa T$, it is possible to extend
$f_\eta=\bigcup\limits_{\alpha<\kappa}f_{\eta\restriction\alpha}$ to an
$\cF$-elementary embedding $f^*$ of $M_{\kappa+1}$ into $N^*$
{\em iff} $\eta\in T$.

Remark that if for example $\chi$ is a strong limit cardinal of cofinality
$\kappa^*$ and $\chi^{<\kappa}\subseteq T\subseteq\chi^{\leq\kappa}\cap
\{\eta^\land\langle0\rangle: (\exists \alpha<\kappa)\lg(\eta)=\alpha+1)\}$,
then
over $\bigcup\limits_{\eta\in\chi^{<\kappa}}M^*_\eta$ for $\chi$ parameters
there are $2^\chi$ independent decisions. This is not only a reasonable
result,  it has been shown, (\cite{Sh:a}, VIII \S1 for $\chi$ as above,
\cite{Sh:e}, III \S5 more generally) that this result is sufficient to prove the
existence of many models in every cardinality $\lambda>\mu+LS(\bT)$.

But to use this construction we have to have some continuity of non forking,
which, we have not proved. Hence we shall use another variant of the
construction

\begin{construction}
\label{const2.2}
{\rm We modify the construction of \ref{const2.1} to suit our purposes.}
\end{construction}

\begin{moddata}
\label{moddata2.2.1}
{\em Let $\langle M_i\in K_{\leq_\mu} : i\leq \kappa+1 \rangle$ be an
continuous $\prenice$-chain of models of $\bT$,
$\|M_{\kappa+1}\|=\mu<\lambda$. Let $T$ be a subset of $
^{\kappa+1\geq}(Ord)$, $<_{lex}$ be the lexicographic order on $T$, it
is a linear order of $T$; suppose
that $T$ is $\triangleleft$-closed i.e. $(\nu\triangleleft\eta\in T\imply
\nu\in T)$, and if $\eta\in {\;}^{\kappa}(Ord)\cap T$, then $\eta^\land\langle
0\rangle$ is the unique $<_{lex}$-successor of $\eta$ in $T$. For $S\subseteq
T$ let $S^{se}=\{\eta\in S:\lg(\eta) \mbox{ successor}\}$. Let $\Op_{i+1}$
wittness
$M_i\prenice M_{i+1}$.  }
\end{moddata}

We define $\Op_\eta=\Op_{\lg\eta}$ for $\eta\in T^{se}$. We can iterate the
operation $\Op_\eta$ w.r.t. $(T^{se}, <_{lex})$.  Also, for each $S\subseteq T$, we
can iterate $\Op_\eta$ w.r.t. $(S^{se}, <_{lex})$.
Let us denote the result of this iteration w.r.t. $(S, <_{lex})$ by $Op^S$
(see \cite[1.11]{KlSh:362}). Note that for any $M\in K$, if
$S_1\subseteq S_2\subseteq
T$, then $M\preceq_{\cF} \Op^{S_1}(M)
\preceq_{\cF}\Op^{S_2}(M)\preceq_{\cF}\Op^T(M)$ (by natural embeddings).
More formally:

\begin{sclaim}
\label{claim2.2.2}
There exist operations $\Op^S$ for $S\subseteq T$ such that
\begin{enumerate}
\item for every $S\subseteq T$ which is $\triangleleft$-closed $M_S=\Op^S(M)$
is defined, and whenever $S_1\subseteq S_2\subseteq T$, then
$M_{S_1}\preceq_{\cF} M_{S_2}$;

\item for $\eta\in T$, $h_{\eta}$ is a surjective $\prec_{\cF}$-elementary
embedding from $M_{\lg(\eta)}$ to $M_\eta$, $M_\eta\preceq_{\cF}
M_{\{\eta\}}$, and $\langle h_{\eta} : \eta \in T \rangle$ is a
$\triangleleft$-increasing 
sequence, i.e. $h_{\eta}\subseteq h_\nu$ whenever $\eta\triangleleft \nu$;

\item for every $x\in M_T$, there exists a $\triangleleft$-closed $S\subseteq
T$, $|S|\leq \kappa$ such that $x\in M_S$ (in fact $S$ is the union of finitely
many branches);

\item for $\eta\in T$ with $\lg(\eta)=\kappa$, and $\alpha<\kappa$, letting
$T[\eta]=\{\nu\in T : \neg(\eta \triangleleft \nu)\}$,
$T^{\leq}[\eta]=\{\nu\in
T[\eta]:\nu\leq_{lex}\eta\}$,
$T^{\geq}[\eta]=\{\nu\in T[\eta]:\eta\leq_{lex}\nu\}$
(so $T[\eta]=T^{\leq}[\eta]\cup T^{\geq}[\eta]$) and
$\nonforkin{M_{T^{\geq}[\eta]}}{M_\eta}_{M_{\eta\restriction\alpha}}^{M_T}$ and
$\nonforkin{M_\eta}{M_{T^{\geq}[\eta]}}_{M_{\eta\restriction\alpha}}^{M_T}$ for
$\alpha<\kappa$;

\item if $\eta\in\lim_\kappa(T)$ and $\eta\notin T$, then
$M_T=\bigcup\limits_{\alpha<\kappa}M_{T[\eta\restriction\alpha]}$

\item $\|M_S\|\leq
|S|+\|M_{\kappa+1}\|^{\kappa^*}+\sup\limits_{\eta\in
S}\|M_{\lg\eta}\|$. 
\end{enumerate}
\end{sclaim}

\begin{sfact}
\label{fact2.2.3}
{\em
\begin{enumerate}
\item By clause (4), if we have the conclusion of \ref{corol1.6} (and
\ref{claim1.18}(1)) then 
$\nonforkin{M_\eta}{M_{T[\eta]}}_{M_{\eta\restriction\alpha}}^{M_T}$.

\item Then in fact one can replace clause (4) above by the weaker condition
\begin{description}
\item[$(4)^-$] for every $S\subseteq T$, if
$\{\eta\restriction i : i\leq \alpha\}\subseteq S \subseteq T$, then
$\nonforkin{M_\eta}{M_S}_{M_{\eta\restriction\alpha}}^{M_T}$.
\end{description}
(2) by (4).
\end{enumerate}
}
\end{sfact}

\noindent{\sc Short Proof of \ref{claim2.2.2}:}\ \ \ \ \
As $\langle M_i:i\leq\kappa+1\rangle$ is $\prenice$-increasing
continuous by renaming
there is $\langle M^*_i:i\leq\kappa+1\rangle$ $\prenice$-increasing
continuous,
$M^*_0=M_0$, $M^*_{i+1}=\Op_{i+1}(M^*_i)$, $M_i\preceq_{\cF} M^*_i$ and
$\nonforkin{M_i^*}{M_{i+1}}_{M_i}^{M^*_{i+1}}$ (for $i\leq \kappa$).
W.l.o.g.
$\|M^*_i\|\leq \|M_i\|^{\kappa^*}$. Let $(I_\eta, D_\eta, G_\eta)$ be
a copy of 
$\Op_\eta$ for $\eta\in T^{se}$ with $I_\eta$'s pairwise disjoint. Define
$I=\Pi\{I_\eta:\eta\in T^{se}\}$, $D,G$ as in the proof of
\cite[1.11]{KlSh:362},
so every equivalence relation $e\in G$ has a finite subset
$w[e]=\{\eta_0^\ell<_{lex}\ldots<_{lex}\eta^\ell_{n(\ell)-1}\}
\subseteq T^{se}$ and ${\frak e}_\ell[e]\in G_{\eta_e^\ell}$ as there. We let
$\Op_{T^{se}}=(I,D,G)$, $M_{T^{se}}=\Op_{T^{se}}(M_0)$ and for $S\subseteq
{\;}^sT$ we let
$$ M_S=\{x\in M_T:w[eq(x)]\subseteq S\}.  $$
Naturally there
are canonical maps $f^*_\eta$ from $M^*_{\lg\eta}$ onto
$M_{\{\nu:\nu\triangleleft \eta\}}$ and let $M_\eta={f''}_\eta(M_{\lg\eta})$.

\begin{improvement}
\label{improv2.2.4}
{\rm Improvement in cardinality.}
\end{improvement}
We can replace $\|M_{\kappa+1}\|^{\kappa^*}$ by $\|M_{\kappa+1}\|+
LS(\bT)$ in part (6) of claim \ref{claim2.2.2}.
After choosing $\langle M^*_i:i\leq\kappa+1\rangle$,
let $M^+_0$ be a Skolemization of $M_0=M^*_0$,
$M^*_{i+1}=\Op(M^+_i)$, $M^+_\delta=\bigcup\limits_{i<\delta}M_i^+$.  Of course
$M^T_S$ ($S\subseteq T$ is $\triangleleft$-closed) are well defined
similarly.  Let $N_i$
be the Skolem hull of $M_i$ in $M^*_i$.  For $\eta\in T$ let
$N_\eta=f^*_\eta(N_{\lg\eta})$.  Now for any
$\triangleleft$-closed $S\subseteq T$
let
$$ N_S=\mbox{Skolem hull in }M^+_S\mbox{ of }\cup\{N_\eta:\eta\in S\}.  $$
\medskip
\centerline{$\ast$\qquad$\ast$\qquad$\ast$}
\medskip

There are two different ways to carry on the construction (under Data
\ref{moddata2.2.1}). We'll consider each in its turn.

\begin{construction}
\label{const2.3}
{\em Recall that it is possible to iterate the
operation $\Op$ with respect to the linear order $(T,<_{lex})$ and this
iteration can be defined as the direct limit of finite approximations. We
shall use different approximations and take the direct limit we obtain the
required operation.

Suppose that $w\subseteq T$ is closed with respect to $\triangleleft$,
(i.e. initial
segment) and is $<_{lex}$-well-ordered. For each approximation $w$ of this
kind, the iterated ultrapower $\Op^w(M_0)$ of $M_0$ with respect to $w$ is
defined as a limit ultrapower and there are natural elementary embeddings into
this limit. The principal difference is that this limit is a little larger
than a limit obtained using only finite approximations. For example, if
$\langle\eta_n:n\leq \omega\rangle$ is a $<_{lex}$-increasing sequence, then in
$\Op^{\eta_\omega}\Bigg(\ldots \Op^{\eta_n}\Big(\ldots\bigg(
\Op^{\eta_0}\big(M_0\big)\bigg)\Big)\Bigg)$, the last operation
$\Op^{\eta_\omega}$
adds elements which are dispersed over all $\Op^{\eta_n}\big(\ldots
\Op^{\eta_0}(M_0)\big)$.
(This is of more interest when the sequence has length $\kappa$.) Now it is
easy to check the symmetry (for $\eta\in {}^{\alpha}\lambda$, $\alpha<\kappa$)
between the $<_{lex}$-successors and $<_{lex}$-predecessors of $\eta$.

We define the embeddings $h_\eta$ for $\eta\in T$ as follows. For
$\eta=\langle\;\rangle$, $h_\eta=\id\restriction M_0$. If
$\eta=\nu^\land\langle i\rangle$, then $\Op^\eta$ acts on
$M_\nu=h_\nu[M_{\lg(\nu)}]$ and we use the commuting diagram:
\bigskip

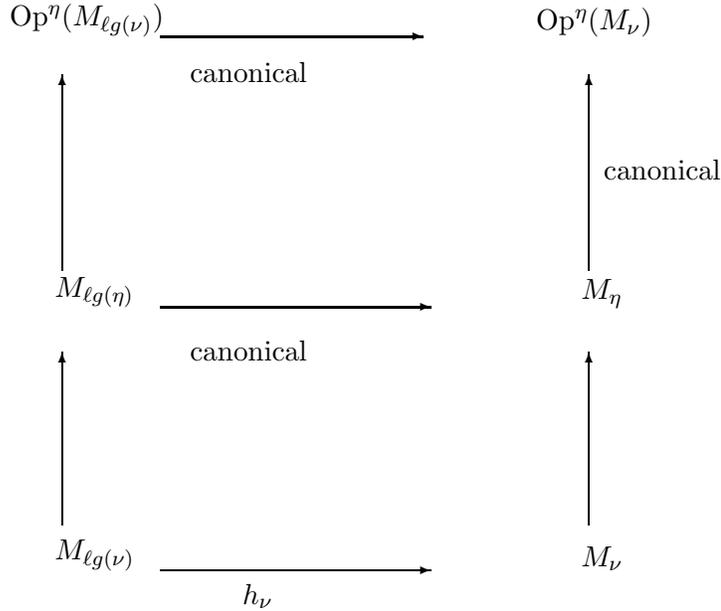
\unitlength=1mm
\special{em:linewidth 0.4pt}
\linethickness{0.4pt}
\begin{picture}(94.00,90.00)
\put(15.00,80.00){\makebox(0,0)[lb]{$\Op^\eta(M_{\lg(\nu)})$}}
\put(85.00,80.00){\makebox(0,0)[lb]{$\Op^\eta(M_\nu)$}}
\put(21.00,44.00){\makebox(0,0)[lb]{$M_{\lg(\eta)}$}}
\put(21.00,9.00){\makebox(0,0)[lb]{$M_{\lg(\nu)}$}}
\put(91.00,44.00){\makebox(0,0)[lb]{$M_\eta$}}
\put(91.00,9.00){\makebox(0,0)[lb]{$M_\nu$}}
\put(35.00,80.00){\vector(1,0){35.00}}
\put(35.00,44.00){\vector(1,0){36.00}}
\put(35.00,9.00){\vector(1,0){36.00}}
\put(22.00,49.00){\vector(0,1){26.00}}
\put(22.00,15.00){\vector(0,1){23.00}}
\put(92.00,49.00){\vector(0,1){26.00}}
\put(92.00,15.00){\vector(0,1){23.00}}
\put(39.00,74.00){\makebox(0,0)[lb]{$\can$}}
\put(39.00,37.00){\makebox(0,0)[lb]{$\can$}}
\put(94.00,61.00){\makebox(0,0)[lb]{$\can$}}
\put(46.00,4.00){\makebox(0,0)[lb]{$h_\nu$}}
\end{picture}

This completes the construction.}
\end{construction}

\begin{construction}
\label{sonstr2.4}
{\em  In this approach, we employ the generalized
Ehrenfeucht-Mostowski models $EM(I,\Phi)$ from chapter VII in
\cite{Sh:a} or \cite{Sh:c}. For this
we need to specify the generators of the model and what the types
are.}
\end{construction}

Let $M_0^+$ be the model obtained from $M_0$ by adding Skolem functions and
individual constants for each element of $M_0$. We know that there is an
operation $\Op$ such that, for $i\leq\kappa$, $M_i\preceq_{\cF}
M_{i+1}\preceq_{\cF} \Op(M_i)$. As in \cite[1.7.4]{KlSh:362} this means that
there are $I,D$ and $G$ such that $\Op(M)=\Op(M,I,D,G)$ where $I$ is a
non-empty 
set, $D$ is an ultrafilter on $I$, and $G$ is a suitable set of equivalence
relations on $I$, i.e.
\begin{description}
\item[(i)] if $e\in G$ and $e'$ is an equivalence relation on $I$ coarser than
$e$, then $e'\in G$;
\item[(ii)] $G$ is closed under finite intersections;
\item[(iii)] if $e\in G$, then $D/e=\{A\subset I/e:\bigcup\limits_{x\in A}x\in
D\}$ is a $\kappa^*$-complete ultrafilter on $I/e$.
\end{description}

For each $b\in M_{i+1}\setminus M_i$, let $\langle x_t^b:t\in I\rangle/D$ be
the image of $b$ in $\Op(M_i)$. We'll also write $\langle x_t^b:t\in
I\rangle/D$ for the canonical image $d(b)$ of $b\in M_i$ in $\Op(M_i)$.

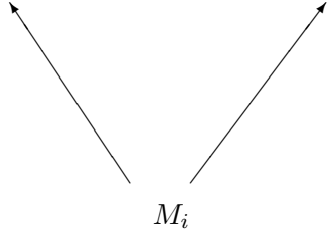
\unitlength=1mm
\special{em:linewidth 0.4pt}
\linethickness{0.4pt}
\begin{picture}(94.00,70.00)
\put(40.00,49.00){\makebox(0,0)[lc]{$M_{i+1}\ni b\ \ \mapsto\ \ \langle x^b_t: t\in I\rangle/D\in \Op(M_i)$}}
\put(71.00,9.00){\makebox(0,0)[lb]{$M_i$}}
\put(68.00,15.00){\vector(-2,3){16.00}}
\put(76.00,15.00){\vector(3,4){18.00}}
\end{picture}

We define a model $M^+$, $M^+_0\preceq_{L_{\kappa^*,\omega}}M^+$, as follows.
$M^+$ is generated by the set $\{x_\eta^b:b\in M_{i+1}\setminus M_i,\; \eta\in
T,\lg(\eta)=i+1\}$. Note that this set does generate a model since $M^+_0$ is
closed under Skolem functions. Since functions have finite arity,
it is enough to
specify, for each finite set of the $x_\eta^b$, what quantifier-free type it
realizes. Since there is monotonicity, we shall obtain indiscernibility as in
\cite{Sh:a}. The type of a finite set $\langle
x_{\eta_\ell}^{b_\ell}:\ell=1,\ldots,n\rangle$ depends on the set $\langle
b_1, \ldots b_n\rangle$
and the atomic (i.e. quantifier-free) type of $\langle
\eta_1,\ldots,\eta_n\rangle$ in the model $\langle
T, \triangleleft ,<_{lex},\mbox{``}\eta\restriction i=\nu\restriction
i\mbox{''}\rangle$. Now w.l.o.g.
we can allow finite sequences $\bar b$ instead of $b$ for $\bar b\in
M_{i+1}\setminus M_i$ and thus w.l.o.g. $\eta_1,\ldots, \eta_n$ is
repetition-free, so w.l.o.g. $\eta_1<_{lex}\eta_2<_{lex}\ldots<_{lex}\eta_n$.
Suppose that the lexicographic order $<_{lex}$ on $\{\eta_\ell\restriction
\alpha: \alpha\leq\lg(\eta_\ell)\mbox{ and }\ell=1, \ldots, n\}$ is a
well-order and the sequence $\langle 
\nu_\zeta:\zeta<\zeta(*)\rangle$ is $\triangleleft$-increasing.  We define
$N_0=M_0^+$, $N_{\zeta+1}=\Op(N_\zeta)$,
$N_\zeta=\bigcup\limits_{\xi<\zeta}N_\xi$ (for limit $\zeta$). Next, we define
$h_{\nu_\zeta}:M_{\lg(\nu_\zeta)}\mathop\to\limits_{\cF} N_{\zeta+1}$,
$h_{\nu_\zeta\restriction\beta}\subseteq h_{\nu_\zeta}$. If $\lg(\nu)$ is a
limit ordinal, then
$\alpha<\lg(\nu)\Rightarrow h_{\nu\restriction\alpha}$
is defined and we
let $h_\nu=\bigcup\limits_{\alpha<\lg(\nu)}h_{\nu\restriction\alpha}$. If
$\nu_\zeta=\nu_\xi\char94\langle\gamma\rangle$, $i=\langle u_\xi)$, then
$M_{\zeta+1}=\Op(M_\zeta,I,D,G)$, identifying elements of $M_\zeta$ with their
images in the ultrapower. Now define
$$
h_{\nu_\zeta}(b)=\left\{\begin{array}{ll} d(H_{\nu_\zeta}(b)) & \mbox{
if } b\in M_i,\\
\langle h_{\nu_\zeta}(x^b_t): t\in I\rangle/D & \mbox{ if } b\in
M_{i+1}\setminus M_i,
\end{array}
\right.
$$
where $d(h_{\nu_\xi}(b))$ is the canonical image of $H_{\nu_\xi}(b)$ in the
ultrapower. The type of $\langle
x_{\eta_\ell}^{b_\ell}:\ell=1,\ldots,n\rangle$ is defined to be the type of
$\langle h_{\eta_\ell}(b_\ell):\ell=1,\ldots,n\rangle$ in $N_\xi$.
\medskip

\begin{sremark}
\label{remark2.4.1}
{\em It is possible to split the construction into two
steps. For $i\leq j\leq\kappa+1$, there is an operation $\Op^{i,j}$,
$M_i\preceq M_j\preceq =\Op^{i,j}(M_i)$, moving $b$ to $\langle
{\;}^{i,j}a_t^b:t\in I\rangle$, $b\in M_j$, ${\;}^{i,j}a_t^b\in M_i$, with the
obvious commutativity and continuity properties. Now the construction is done
on a finite tree $\langle\eta_\ell:\ell=1,\ldots,n\rangle$, $\langle
\eta_\ell\cap\eta_m:\ell,m<\omega\rangle$.  We omit the details of
monotonicity.}
\end{sremark}

\begin{snotation}
\label{notation2.4.2}
{\em Let $M_T=M$ be the Skolem closure. If $S\subseteq
T$ is closed with respect to initial segments, let
$M_S=Sk_{M_T}(x_\eta^b:\eta\in S,b\in M_{\lg(\eta)})$ and
$M^*_\eta=M_{\{\eta\restriction\alpha: \alpha\leq\lg(\eta)\}}$. Define
$h_\eta:M_{\lg(\eta)}\to M^*_{\eta}$ by
$h_\eta(b)=x^b_{\eta\restriction\tau(\bT)}$ and
$N_\eta=h_\eta[M_\eta]$.}
\end{snotation}

\begin{ssremark}
\label{remark2.4.3}
{\em The construction can be used to get many fairly
saturated models. We list the principal properties below.}
\end{ssremark}

\begin{ssfact}
\label{fact2.4.4}
Suppose that $S_\ell\subseteq T$ is closed with respect to initial segments,
$S_0=S_1\cap S_2$ and 
$[\eta\in S_1\ \&\ \nu\in S_2\setminus S_1\imply
\eta<_{lex}\nu]$  then 
$$
\nonforkin{M_{S_1}}{M_{S_2}}_{M_{S_0}}^{M_T}.
$$
\end{ssfact}

\Proof
W.l.o.g. $S_\ell$ is closed, $M_{cl(S_\ell)}=M_{S_\ell}$. Let
$S_2\setminus S_0=\{\nu_\zeta:\zeta<\zeta(*)\}$ be a list such that
$\nu_\zeta<\zeta_\xi\imply\zeta<\xi$; let
$S_2^\xi=S_0\cup\{\nu_\zeta:\zeta<\zeta(*)\}$. Then
\begin{enumerate}
\item $\langle M_{S_2^\xi}:\xi\leq\xi(*)\rangle$ is continuous increasing;
\item $\langle M_{S_2^\xi\cap S_1}:\xi\leq\xi(*)\rangle$ is continuous
increasing.
Hence one has
\item $\nonforkin{M_{S_2^\xi\cup
S_1}}{M_{S_2^{\xi+1}}}_{M_{S_2^\xi}}^{M_{S_2^{\xi+1}\cup S_1}}$
\end{enumerate}

This is immediate from the definitions, because $M_{S_2^{\xi+1}\cup S_1}$ is
the Skolem closure of $M_{S_\xi^2\cup S_1}\cup N_{\nu_\xi}$, and so elements of
$N_{\nu_\xi}$ can be represented as averages.

\section{Categoricity in $\mu$, when $LS(\bT)\leq \mu < \lambda$}

\begin{hypothesis}
\label{hyp3.0}
Every $M\in K_{<\lambda}$ is nice hence has a $\prec_{\cF}$-extension
of cardinality $\lambda$ which is saturated and ${\cal K}_{<\lambda}$
has amalgamation.
\end{hypothesis}

This section contains the principal theorems of the paper: if $\bT$ is
$\lambda$-categorical, $LS(\bT)\leq\mu<\lambda$, then
$\boldk_\mu(\bT)=\emptyset$ when $\mu\in [LS(\bT),\lambda)$ and
when $LS(\bT)\leq \chi=\cf\chi<\lambda$, $\bT$ is $\chi$-based, (and
${\cal K}$ does not have 
$(\mu,\kappa)$-continuous non forking when $\mu\in[LS(\bT),\lambda]$,
$\kappa\leq\mu$) also there is a saturated model in
${\cK}_\mu=\langle K_\mu,\preceq_{\cF}\rangle$ and $\bT$ is
$\lambda$-categorical. However we first deal with some preliminary
results, quoting \cite{Sh:300} extensively.

\begin{theorem}
\label{th3.1}
Assume the conclusion of \ref{corol1.6} for $\mu$ (e.g. $\mu^+<\lambda$).
Suppose that the tree $T$ is as in Claim \ref{claim2.2.2} and
 suppose
further: $\langle M_i\in K_{\leq\mu}:i\leq\kappa+1\rangle$ is
$\prenice$-increasing continuous sequence of members of $K_{\leq
\mu}$, and we apply \S2 and 
\begin{description}
\item[$(*)$] there is no $\preceq_{{\cal F}}$-increasing continuous
sequence $\langle N_i\in K_{\leq\mu}:i\leq\kappa\rangle$ such
that:
$$ M_i\preceq_{\cal F} N_i$$
$$ M_{\kappa+1}\preceq_{\cal F} N_\kappa$$
$$\nonforkin{N_i}{M_{i+1}}_{M_i}^{N_{i+1}}\hbox{ for }i<\kappa$$
\end{description}
{\em Then} TFAE for $\eta\in \Lim_\kappa(T)\stackrel{\rm
def}{=}\{\eta\in {\;}^{\kappa}(Ord):  
\bigwedge\limits_{i<\kappa}(\eta\restriction (i+1)\in T)\}$:
\begin{description}
\item[$(\alpha)$] There is an ${\cF}$-elementary embedding $h$ from
$M_{\kappa+1}$ into $M_T$ such that
$\bigcup\limits_{i<\kappa}h_{\eta\restriction i+1}\subseteq h$.
\item[$(\beta)$] $\eta\char94\langle 0\rangle \in T$.
\end{description}
\end{theorem}

\Proof As regards the implication from $(\beta)$ to $(\alpha)$, so assume
$\eta\in T$ consider the ${\cal F}$-elementary embedding $h_{\eta^\land\langle
0\rangle}$. Check that $h_{\eta^\land \langle 0\rangle}$ is as required in
$(\alpha)$. The other direction follows by \ref{fact2.2.3}(1) and $(*)$.
\hfill\qed$_{\ref{th3.1}}$

\begin{claim}
\label{claim3.1A}
Suppose the conclusion of \ref{corol1.6} for $\mu$ and
$\bar M=\langle M_i\in K_{\leq\mu}:i\leq\kappa+1\rangle$ is given.
Then $\bar M$ satisfies $(*)$ of
\ref{th3.1} if one of the following holds:
\begin{description}
\item[$(\alpha)$] there is $a\in M_{\kappa+1}$ such that
$i<\kappa \imply
M_\kappa\doesforkin_{M_{i+1}}^{M_{\kappa+1}} a$, or
\item[$(\beta)$] $\kappa=\cf(\kappa)=\mu > LS(\bT)$
and $i<\kappa \imply \|M_i\|<\kappa$, and
there is a continuous $\prec_{\cF}$-chain $\langle N_i :
i\leq \kappa\rangle$,  $M_{\kappa+1} = \bigcup\limits_{i\leq \kappa}N_i$,
$\kappa=\chi^{\cf(\kappa)}$, $\bigwedge\limits_{i<\kappa}(N_i\in
K_{<\kappa})$, and $E=\{i<\kappa: M_{i+1}\doesforkin_{M_i}^{N_\kappa}
N_i\}$ is a stationary subset of $\kappa$.
\end{description}
\end{claim}

\Proof Straight from \ref{th3.1}, \S2.
\medskip

\begin{remark}
\label{rem3.1B}
{\em Clause $(\beta)$ can also be proved using niceness as in the proof of
\ref{th3.4}. This works for any $\kappa<\lambda$. Also we can 
imitate \ref{claim2.2.2} but no need arise.}
\end{remark}

\begin{corollary}
\label{corol3.2}
If $\bT$ is a $\lambda$-categorical theory\footnote{or just has $<
2^\lambda$ non isomorphic model in $\lambda$}, {\em then}
\begin{enumerate}
\item $\bT$ is $\chi$-based if $\chi^+<\lambda$ and $\chi\geq
LS(\bT)$; also it is not $(<\mu)$-based if
$\mu=\cf(\mu)$, $LS(T)<\mu$, 
$\mu^+<\lambda$;
\item $\boldk_\mu(${\bf T}$)=\emptyset$ for every $\mu$, $\mu^+<\lambda$
and $\mu\geq LS(\bT)$.
\end{enumerate}
\end{corollary}

\Proof
1), 2) We use \ref{th3.1}, \ref{claim3.1A} to contradict
$\lambda$-categoricity. 

\noindent{\sc Case 1}:\ \ \ $\lambda^\mu=\lambda$ By \cite{Sh:300},
III, 5.1 = \cite{Sh:e} IV, 2.1.
\medskip

\noindent{\sc Case 2}:\ \ \ $\lambda$ is regular, $\lambda>\mu^+$.
We can find a stationary $W^*\in I[\lambda]$, $W^*\subseteq
\{\delta<\lambda: \cf(\delta)=\kappa\}$ (by \cite{Sh:420}, \S1). Hence,
possibly replacing $W^*$ by its intersection with some club of
$\lambda$, there is $W^+$, $W^*\subseteq W^{++}$ and $\langle a_\alpha:
\alpha\in W^+\rangle$ such that: $\alpha\in a_\beta$ (so $\beta\in W^+$)
implies $\alpha\in W^+$, $a_\alpha=a_\beta\cap a_\alpha$ and $\otp
(a_\alpha)\leq \kappa$ and $\alpha=\sup a_\alpha \iff
\cf(\alpha)=\kappa \iff \alpha\in W^*$.
Now let $\eta_\alpha$ enumerate $a_\alpha$ in increasing order (for
$\alpha\in W^+$), and for any $W\subseteq W^*$ let
$$
T_W=\{\eta_\alpha:
\alpha\in W^+ \mbox{ but } \alpha\notin W^*\setminus W\}\cup
\{\eta_\alpha\char94\langle 0\rangle: \alpha\in W\}.
$$
Now if $W_1$, $W_2\subseteq W$, $W_1\setminus W_2$ is stationary, then
$M_{T_{W_1}}$ cannot be $\preceq_{\cF}$-embedded into
$M_{T_{W_2}}$(again by \cite{Sh:300} III, \S5 = \cite{Sh:e}, IV \S2).
\medskip

\noindent{\sc Case 3}:\ \ \ $\lambda$ singular.

Choose $\lambda'$, $\lambda>\lambda'=\cf(\lambda')>\mu^+$ and act as
in case 2 (to get $2^\lambda$ we need more, see \cite{Sh:e}, IV.
\hfill\qed$_{\ref{corol3.2}}$

\begin{hypothesis}
\label{hyp3.5A}
The conclusion of \ref{corol3.2} (in addition to \ref{hyp3.0} of course).
\end{hypothesis}

\begin{conclusion}
\label{conc3.3}
Suppose $\mu\geq LS(\bT)$, $\mu^+<\lambda$, $M\in K_\mu$
\begin{enumerate}
\item If $p\in S(M)$ then $p$ is determined by $\{p\restriction
N:N\preceq_{\cF} M\mbox{ and } \|N\|=LS(\bT)\}$

\item Assume further
\begin{description}
\item[$(*)^M_{\{N_t:t\in I\}}$\ a)]
 $I$ (a partial order) which is directed (i.e. every finite 

$\qquad\quad$  many elements have a common upper bound)
\item[\qquad \qquad b)] $N_t\preceq_{\cF} M$,
\item[\qquad \qquad c)] $I\models t\leq s$ implies $N_t\subseteq N_s$
\item[\qquad \qquad \quad\ ] (hence $N_t\preceq_{{\cal F}} N_s$ by
clause (b)) 
\item[\qquad \qquad d)] $\bigcup\limits_{t\in I}M_t=M$.
\end{description}
\end{enumerate}
{\em Then} every $p\in S(M)$ is determined by
$$
\{p\restriction N_t:t\in I\}
$$
\end{conclusion}

\Proof 1) Follows by part (2).

\noindent 2) Easily (and as \cite{Sh:88} \S1):
\begin{description}
\item[$\otimes$] we can choose by induction on $n<\omega$ for every
$u\in[M]^n$,
$t[u]\in I$ and $N^*_{u}$ such that:
\begin{description}
\item[\ ] $u\in N^*_{u}$, $N^*_{u}\preceq_{\cF} N_{u}$,
$\|N^*_{u}\|\leq LS(\bT)$ and:
$u\subseteq v\in [|M|]^{<\aleph_0}$ implies $N^*_{u}\prec N^*_{v}$ and
$t[u]\leq_I t[v]$.
\end{description}
\end{description}

Let for $U\subseteq |M|$, $N^*_{U}=:\cup\{N^*_u:u\subseteq U\mbox{
finite}\}$ the definitions are compatible.
Easily $U_1\subseteq U_2\subseteq |M|\imply N^*_{U_1}\preceq_{\cF}
N^*_{U_2}\preceq_{\cF} M$.  Now we prove by induction on $\mu\leq\|M\|$ that:
\begin{description}
\item[$(**)$] if $U\subseteq \|M\|$, $|U|=\mu$, $p\in N^*_U$ then for some
$u\in [U]^{<\aleph_0}$, $p$ does not fork over $N^*_u$.
\end{description}

For $\mu$ finite this is trivial, for $\mu$ infinite then $\cf(\mu)\notin
\boldk_{\mu+LS(\bT)}(\bT)$ (by \ref{corol3.2}(2)) so $(**)$ holds. Now by
\ref{claim1.20}(3), we are done.
\hfill\qed$_{\ref{conc3.3}}$

\begin{theorem}
\label{th3.4}
Suppose that $\cf(\kappa)=\kappa\leq\mu<\lambda$ 
 and $LS(\bT) < \mu$.  Then
\begin{enumerate}
\item The $(\mu,\kappa)$-saturated model $M$ is saturated (i.e. 
$N\preceq_{\cF} M$, $\|N\|< \|M\|$, $p\in S(N)\imply p$ realized in
$M$,  
and hence unique). Hence there is a saturated model in $K_\mu$.
\item The union of a continuous $\preceq_{\cF}$-chain of length $\kappa$ of
saturated models from $K_\mu$ is saturated.
\item In part (1) we can replace saturated by $(\mu, \mu)$-saturated
if $\mu= LS(\bT)$.\\ 
We can in part (1) replace saturated by $\chi$-saturated
if $\chi> LS(\bT)$ .

\end{enumerate}
\end{theorem}

\Proof 1), 2) Suppose that $M=M_*$ and $\langle M_i:i\leq\kappa\rangle$ is a
continuous $\preceq_{\cF}$-chain of members of $K_\mu$
such that for the proof of 1) $M_{i+1}$ is a
universal extension of $M_i$ and for the proof of 2) $M_{i+1}$ is saturated. 
Let $i\leq j\leq\kappa$. Then $M_i\prenice M_j$
(by \cite{KlSh:362}, 5.4 or more exactly by the hypothesis \ref{hyp3.0}). So there is an operation $\Op_{i,j}$ such that
$M_i\preceq_{\cF} M_j\preceq_{\cF} \Op_{i,j}(M_i)$. It follows that there
is an expansion $M^+_{i,j}$ of $M_j$ by at most $LS(\bT)$ Skolem functions
such that if $N$ is a submodel of $M^+_{i,j}$, then
$$
\nonforkin{M_i}{M_j\restriction N}_{M_j\restriction(N\cap M_i)}^{M_j}.
$$
[Why? as we use operations coming from equivalence relations with $\leq
\kappa^*$ classes and $LS(\bT)\geq \kappa^*$ by its definition]. More fully,
letting $\Op_{i, j}(N)= N^I_D/G$, 
every element $b\in M_j$ being in $\Op_{i, j}(M_i)$ has a
representation as the equivalence class of $\langle x_t^b:t\in I\rangle/D$
under $\Op_{i,j}$, $x_t^b\in M_i$ and $|\{x_t^b:t\in I\}|\leq\kappa^*$. The
functions of $M^+_{i,j}$ are the Skolem functions of $M_j$ and $M_i$
and functions $F_\zeta$ ($\zeta<\kappa^*$) such that $\{F_\zeta(b):
\zeta< \kappa^*\}\supseteq \{x^b_t: t\in I\}$.

If $\kappa=\cf(\mu)$, the theorem is immediate. So we'll suppose that $\kappa<\mu$.
Suppose $N\preceq M=M_\kappa$, $\|N\|<\mu$ and $p\in S(N)$. Let
$\chi=: \|N\|+\kappa+LS(\bT)$.
W.l.o.g. there is no $N_1$, $N\preceq_{\cF} N_1\prec M_\kappa$,
$\|N_1\|\leq \chi$ and $p_1$, $p\subseteq p_1\in S(N_1)$ such that $p$ fork
over $N$ (by \ref{claim3.1A}). If there is $i <
\kappa$ such that $N\subseteq M_i$, then $p$ is realized in $M_{i+1}$.
By the choice of the models $M^+_{i,j}$, it is
easy to find $N'$ such that $N\preceq N'\preceq M_\kappa$,
$\|N'\|=\chi\stackrel{\rm def}{=}\|N\|+\kappa+LS(\bT)$ and,
for every $i\leq\kappa$,
$$
\nonforkin{M_i}{N'}_{M_i\cap N'}^{M_\kappa}.
$$
Now let $N_i=N'\cap M_i$ and note that $N_\kappa=N'$. The sequence $\langle
N_i:i\leq\kappa\rangle$ is continuous increasing and there is an extension
$p'$ of $p$ in $S(N_\kappa)=S(N')$. Hence there exists $i<\kappa$ such that
$(i\leq j<\kappa)\imply$ ($p'$ does not fork over $N_j$). By
\ref{corol3.2}(1), it is
sufficient to find $j\in [i,\kappa)$
and $\langle N^*_\varepsilon: \varepsilon<\chi^+\rangle$
$\preceq_{\cF}$-increasing continuous such that: 
$N_i\preceq_{\cF} N^*_\varepsilon \preceq_{\cF} M_j$, 
$N^*_{\varepsilon+1}$ is a $\chi$-universl extension
of $N^*_\varepsilon$
(recall symmetry and uniqueness of extensions). 

3) Similar proof for the second sentence, \ref{conc1.17} for the first
sentence.
\hfill\qed$_{\ref{th3.4}}$
\medskip

\noindent
{\bf Remark:} Using categoricity we can prove \ref{th3.4} also by
\ref{conc1.17}(2) (and uniqueness).

\medskip
\begin{conclusion}
\label{conc3.5}
Assume $LS(\bT)\leq\kappa<\mu\in(LS(\bT),\lambda)$, $M\in
K_\mu$ is not $\kappa^+$-saturated; let $\langle N^*_u:u\in
[|M|]^{<\aleph_0}\rangle$ and $N^*_U$ (for $U\subseteq |M|$) be as in
the proof of \ref{conc3.3}(2).
{\em Then} there is $U\subseteq |M|$, $|U|\leq\kappa$, $p\in S(N^*_U)$ i.e.
there are $N^+$, $N^*_U\preceq_{\cF} N^+\in K_\kappa$, and $a^+\in
N^+$ satisfying $(a^+, N^+)/E_{N^*_U}=p$ such that:\\
 for no $a\in M$ do
we have $u\in 
[U]^{<\aleph_0}\imply \tp(a,N^*_u,M)=\tp(a^+,N^*_u,N^+)$. Equivalently:
w.l.o.g. $N^+\cap M=N^*_U$ and we can define $N^+_u$ for
$u\in [|N^+|]^{<\aleph_0}$,
such that $\langle N^+_u:u\in [|N^+|]^{<\aleph_0}\rangle$
as in the proof of \ref{conc3.3}(2), 
and $u\in [U]^{<\aleph_0}\imply N^+_u=N^*_u$
and for no $u_0\in [|M|]^{<\aleph_0}$, $v_0\in
[|N^+|]^{<\aleph_0}$, $a^+\in N^*_{v_0}$,
and $a\in N_{u_0}^*$ do we
have
$$
\bigwedge\limits_{u\in [U]^{<\aleph_0}}\tp(a,N^*_u,N^*_{u\cup u_0})
=\tp(a^+, N^+_u,N^+_{u\cup v_0}).
$$
\end{conclusion}

\begin{corollary}
\label{corol3.6}
\begin{enumerate}
\item If $\bT$ is $\lambda$-categorical and $LS(\bT)< \mu< \lambda$,
$LS(\bT)\leq\chi$, $\delta(*)=(2^{LS(\bT)})^+$ and
$\beth_{\delta(*)}(\chi)\leq\mu$ {\em then} every
$M\in K_\mu$ is $\chi^+$-saturated. In fact for some
$\delta<\delta(*)$ we can replace $\delta(*)$ by $\delta$.
\item If $\mu=\beth_{(2^\chi)^+\times \delta}$,
$\delta$ a limit ordinal {\em then}
$\bT$ is $\mu$-categorical.
\end{enumerate}
\end{corollary}

\Proof By \ref{conc3.5} this problem is translated to an omitting
type argument + cardinality of a predicts which holds
(see \cite{Sh:c}, VIII \S4, \cite{Sh:c}, VII \S5). See more on this in
\cite{Sh:88}.
\hfill\qed$_{\ref{corol3.6}}$
\medskip

\begin{claim}
\label{claim3.7}
[$\bT$ categorical in $\lambda$]
\begin{enumerate}
\item If $\langle M_i: i\leq \delta\rangle$ is
$\preceq_{\cF}$-increasing continuous, $M_i\in K_{<\lambda}$, $p\in
S(M_\delta)$ {\em then} for some $i<\delta$, $p$ does not fork over $M_i$.
\item If $N\in K_{<\lambda}$ and $p$, $q\in S(N)$ does not fork over
$M$, $M\preceq_{\cF} N\in K_{<\lambda}$ {\em then} $p=q\iff
p\restriction M=q\restriction M$. Moreover if
$M\preceq_{\cF}N\preceq_{\cF} N^+$, $a\in N^+$ then 
$$
\nonforkin{N}{a}_{M}^{N^+} \Leftrightarrow \nonforkin{a}{N}_{M}^{N^+}.
$$
\item If $M\preceq_{\cF} N\in K_{<\lambda}$ and $p\in S(M)$ {\em then}
there is $q\in S(N)$ extending $p$ not forking over $M$.
\item If $M_0\preceq_{\cF} M_1\preceq_{\cF} M_2$, $p\in S(M_2)$,
$p\restriction M_{\ell+1}$ does not fork over $M_\ell$ for $\ell=0,1$
{\em then} $p$ does not fork over $M_0$.
\item If $\mu$, $\delta<\lambda$, $M_i\in K_{\leq \mu}$ for $i<\delta$
is 
$\preceq_{\cF}$-increasing continuous, $p_i\in S(M_i)$, $[j<i \Rightarrow
p_j\subseteq p_i]$ {\em then} there is $p\in S(M_\delta)$ such that $i<\delta
\Rightarrow p_i\subseteq p_\delta$.
\end{enumerate}
\end{claim}

\Proof
1) Otherwise we can find $N$, $M_\delta\preceq_{\cF} N\preceq_{\cF}
\Op(M_\delta)$, $N\in K_{\lambda}$, $N$ omit $p$:
$\bigcup\limits_{i<\delta}\Op(M_i)$; so we get a non $\lambda$-saturated
model of cardinality $\geq \lambda$, contradiction.

2) The first sentence follows from the second. If the second fails
then we can 
contradict stability in $\|N\|$, by a proof just like
\ref{th1.5}, \ref{corol1.6}.

3) we can find an operation $\Op$, $\|\Op(M)\|\geq \lambda$, so in
$\Op(M)$ some $\bar a$ realizes $p$ so $q=\tp(\bar a, N, \Op(N))$ is as
required.

4) For some operation $\Op$, some $\bar a\in {}^{\omega>}(\Op(M_0))$
realizes $p\restriction M_0$, so $p_\ell=\tp(\bar a, M_\ell, \Op(M_l))$
does not fork over $M_0$, and $p_{\ell+1}$ does not fork over $M_\ell$, so by part 2) show
$p_1= p\restriction M_1$ and then $p_2=p$, but $p_2$ does not fork over $M_0$,

5) {\em Case 1:} $\cf(\delta)>\aleph_0$

\noindent
For every limit $\alpha<\delta$ for some $i<\delta$ we have $p_\delta$
does not fork over $M_\alpha$. By Fodour lemma,
for some $i<\delta$, $j\in [i, \delta) \Rightarrow p_j$ does not fork over
$M_i$. So the stationarization of $p_i$ in $S(M_\delta)$ (exists by
\ref{claim1.19}) is as required.

\noindent
{\em Case 2:} $\cf(\delta)=\aleph_0$.

\noindent
So w.l.o.g. $\delta=\omega$. Here chasing arrows (using amalgamation)
suffice.\\
\   
\hfill\qed$_{\ref{claim3.7}}$
\medskip

\begin{lemma}
\label{lemma3.8}
In $K_{<\lambda}$ we can define
$\rrk(\tp(a,M,N))$ with the right properties. I.e.
\begin{description}
\item[(A)] if $M\prec_{\cF} N\in K_{<\lambda}$, $\bar a\subseteq N$,
$M\in \bigcup\limits_{\mu^+<\lambda}K_\mu$, $p=\tp(\bar a, M, N)$ then
$$
\begin{array}{ll}
\rrk(p)\geq \alpha \mbox{ iff } & \mbox{ for every }\beta<\alpha\mbox{
there are }\\
\ & p', M' \mbox{ such that } M\prec_{\cF} M'\in
\bigcup_{\mu^+<\lambda}K_\mu\\
\ & p'\in S(M'), p'\restriction M=p\mbox{ and }\rrk(p')\geq \beta
\end{array}
$$
\item[(B)] for every $M$, $N$, $\bar a$, $p$ as above $\rrk(p)$ is an
ordinal.
\item[(C)] If $M_1\prec_{\cF} M_2\in \bigcup\limits_{\mu^+<\lambda}
K_\mu$, $p_2\in S(M_2)$, then $\rrk(p_2\restriction M_1)\geq
\rrk(p_2)$ and equality holds iff $p_2$ does not fork over $M_1$ and
then $p_2\restriction M_1$ (and $M_2$) determine $p_2$
\item[(D)] If $\langle M_i: i\leq \delta\rangle$ is
$\preceq_{\cF}$-increasing continuous, $M_i\in
\bigcup\limits_{\mu^+<\lambda} K_\mu$ and $p_\delta\in S(M_\delta)$
then for some $i<\delta$ we have: $j\in [i, \delta]\Rightarrow
\rrk(p_\delta)=\rrk(p_\delta\restriction M_j)$.
\end{description}
\end{lemma}

\begin{lemma}
\label{lemma3.9}
Assume $\mu\geq LS(\bT)$, $\mu^+<\lambda$.  If
$M\in K_\mu$ is saturated (for $\mu=LS(\bT)$ means
$(\mu,\mu)$-saturated), and
$p\in S(M)$ {\em then} there are $N$, $a$ such that $N\in K_\mu$ is
saturated, $a\in N$, 
$\tp(a,M,N)=p$ and $N$ is $\mu$-isolated over $M\cup\{a\}$ (i.e. if
$N\preceq_{\cal F} N^+\in K_{<\lambda}$ and $N^*\preceq_{\cF}
N^+$,
and $\tp( a, N^*, N^+)$ does not fork over $N$ ($\preceq_{\cF} N^*$) {\em then}
$\nonforkin{N^*}{N}_{M}^{N^+}$).
\end{lemma}

\Proof  As in \cite{Sh:h}, Ch. V (or Makkai Shelah \cite{MaSh:285},
4.22) because we have \ref{corol3.2}(1) (by \ref{hyp3.5A}).

\hfill\qed$_{\ref{lemma3.9}}$
\medskip

\begin{claim}
\label{claim3.10}
For $M$, $a$, $N$ as in \ref{lemma3.9}, if $N\preceq_{\cF} N^+\in
K_{<\lambda}$, $A\subseteq N^+$, $|A|\leq \mu$ and
$\nonforkin{a}{A}_{M}^{N^+}$ {\em then} $\nonforkin{N}{A}_{M}^{N^+}$.
\end{claim}

\Proof We use the symmetry of $\nonforkin{}{}_{}^{}$ (hold by
\ref{corol1.6} as $\mu^+<\lambda$).
\hfill\qed$_{\ref{claim3.10}}$
\medskip

\begin{claim}
\label{claim3.11}
If $\mu\in [LS(\bT), \lambda)$, $M\in K_\mu$ is saturated and $p\in
S(M)$ {\em then} for some saturated $N\in K_\mu$, $M\preceq_{\cF} N$, $a\in N$
and $(M, N, a)$ satisfies the conclusion of \ref{claim3.10} for finite
$A$.
\end{claim}

\Proof A problem arise only if $\mu^+=\lambda$. We can find
$\langle M'_i: i\leq \mu\rangle$ which is $\preceq_{\cF}$-increasing
continuous, $\|M'_i\|=|i|+ LS(\bT)$, $M'_\mu=M$, $M'_i$ is saturated,
$M'_{i+1}$ universal over $M'_i$ and $p$ does not fork over $M_0$.

Now choose by induction on $i\leq \mu$, $(M_i, N_i, a)$ such that:
\begin{description}
\item[(a)] $M_0=M'_0$,
\item[(b)] $\|M'_i\|=\|N'_i\|=|i|+ LS(\bT)$,
\item[(c)] for $i$ non limit $(M_i, N_i, a)$ as in \ref{lemma3.9} (with
$|i|+LS(\bT)$ instead $\mu$),
\item[(d)] $\tp(a, M_0, N_0)=p\restriction M'_0$,
\item[(e)] $\langle M_i: i\leq \mu\rangle$ is $\preceq_{\cF}$-increasing
continuous,
\item[(f)] $\langle N_i: i\leq \mu\rangle$ is $\preceq_{\cF}$-increasing
continuous,
\item[(g)] $\tp(a, M_{i+1}, N_{i+1})$ does not fork over $M_i$ (hence is the
stationarization of $\tp(a, M_0, N_0)=p\restriction M'_0$),
\item[(h)] $M_{i+1}$ is universal over $M_i$.
\item[(i)] $M_i\preceq_{\cF} N_i$.
\item[(j)] $N_{i+1}$ is isolated over $(M_{i+1}, a)$
\end{description}
There is no problem, so as $M_\mu$ is saturated and in $K_\mu$,
$M_0 = M'_0$ has cardinality $<\mu$, w.l.o.g. $M_\mu=M$. For any
candidates $N^+$, $A$, as in \ref{claim3.10} (but $A$ is finite) assume
$N\doesforkin_{M}^{N^+}A$; as $A$ is finite, for some $i<\mu$, the
type $\tp(A, M, N^+)$ does not fork over $M_i$, and for some $j<\mu$ the type
$\tp(A, N, N^+)$ does not fork over $N_j$, w.l.o.g. $i=j$ is a
successor ordinal and $\tp(A\cup\{a\}, M)$ does not fork over $M_i$.
So as $N\doesforkin_{M}^{N^+} A$, neccessarily $\tp(A, N_i, N^+)$
forks over $M_i$, hence (by clause (c) above),
$a\doesforkin_{M_i}^{N^+}A$, hence $a \doesforkin_{M}^{N^+}A$
(state the laws of $\nonforkin{}{}_{}^{}$).

Alternatively repeat the proof of \ref{lemma3.9} using
\ref{claim3.7}(2)'s second sentence.
\hfill\qed$_{\ref{claim3.11}}$
\medskip

\begin{theorem}
\label{th3.12}
Assume $\lambda$ is a successor cardinal i.e. $\lambda=\lambda_0^+$.  Then
$\bT$ is categorical in every $\mu\in [\beth_{(2^{LS(\bT)})^+},\lambda)$
(really for some $\mu_0<\beth_{(2^{LS(\bT)})^+}$, $\mu\in
[\mu_0,\lambda)$ suffices).
\end{theorem}

\Proof As in \cite{MaSh:285}.
By \ref{corol3.6}, for some $\mu_1<\beth_{(2^{LS(\bT)})^+}$
every $M\in K_{[\mu_1,\lambda]}$ is
$LS(\bT)^+$-saturated. Let $\mu\in [\mu_1,\lambda)$, and assume
$M\in K_\mu$ is not saturated, so for some $\kappa\in (LS(\bT),\mu)$ the
model
$M$ is $\kappa$-saturated not $\kappa^+$-saturated. Let $p$, $\langle
N^*_u:u\in [|M|]^{<\aleph_0}\rangle$, $U$, $N^+$, $\langle N^+_u:u\in
[|N^+|]^{\aleph_0}\rangle$ be as in \ref{conc3.5}. Let $U_0=U$.
W.l.o.g. $N^*_{U_0}$ is saturated, $p$ does not fork over $N_{u^*}^*$,
$u^*\in [U]^{<\aleph_0}$
finite, $\rrk(p)$ minimal under the circumstances. Now let $b\in M\setminus
N^*_{U_0}$, so there is $N_1\preceq_{\cal F} M$ which is $\mu$-isolated over
$N^*_{U_0}\cup\{b\}$. By defining more  $N^*_u$ w.l.o.g.
$N_1=N^*_{U_1}$. So $\tp(b, N^*_{U_0}, M)$, and $p$ are orthogonal (see
\cite{Sh:h}, Ch. V).
Now we deal with orthogonal types and we continue as \cite{MaSh:285}:
define a $\prec_{\cF}$-chain $M^*_i$ ($i<\lambda$) of saturated models
of cardinality $\lambda_0$ all omitting some fixed $p\in S(M^*_0)$.
\hfill\qed$_{\ref{th1.13}}$
\medskip

\begin{discussion}
\label{3.15}
{\rm 
1) Below $\beth_{(2^{LS(\bT)})^+}$
}
\end{discussion}

A problem is what occurs in $[LS({\bf T}),\beth_{(2^{LS(\bT)})^+}]$.
As $\bT$ is not necessarily complete, for any $\psi$ and $\bT$ we can consider
${\bT}'\stackrel{\rm def}{=}\{\psi \rightarrow \varphi:\varphi\in \bT\}$,
if
$\neg\psi$ has a model in $\mu$ iff $\mu<\mu^*$, we get such examples.  So we
may consider $\bT$ complete.  Hart Shelah \cite{HaSh:323} bound our possible
improvement but we may want larger gaps, a worthwhile direction.  If
$\bT\subseteq L_{\kappa^+,\omega}$ is $L_{\kappa^+,\omega}$-complete hence
$L_{\infty,\omega}$-complete, $LS(\bT)=\kappa$,  we cannot improve.

If $|\bT|<\kappa^*$ we may look at what occurs in large enough
$\mu<\kappa^*$.
\medskip

\noindent 2) Below $\lambda$.
\medskip

If $\lambda$ is a limit cardinal we get only \ref{claim3.7}, this
is a more serious issue.
The problem is that we can get $\mu$-saturated not saturated model in
$K_{\mu^+}$, so we get for $M\in K_\mu$ saturated, two orthogonal types
$p$, $q\in S(M)$ (not realized in $M$). We want to build a prime model over
$M\cup$(a large indiscernible set for $p$). Clearly ${\cal P}^-(n)$-diagrams
are called for.
\medskip

\noindent 3) Above $\lambda$
\medskip

In some sense we
know
every model is saturated: if $M\in K_{>\lambda}$, $N\preceq_{\cal F} M$,
$\|N\|<\lambda$, $p\in S(N)$ then $dim(p,N,M)=\|M\|$,
i.e. if $N\preceq_{\cF} N^+\preceq_{\cF} M$ and : $\|N^+\|<\|M\|$
when $\lambda$ is successor,
or $\beth_{(2^{LS(\bT)})^+}(\|N^+\|)$ when $\lambda$ is a limit cardinal.

Another way to say it: the stationarization of $p$ over $N^+$ is realized. But
is every $q\in S(N^+)$ a stationarization of some $p\in S(N')$,
$N'\preceq_{\cF} N^+$, $\|N'\|\leq LS(\bT)$?  We can find
$N_0\preceq_{\cF} N^+$, $\|N_0\|\subseteq (\bT)$, such that:
$[N_0\preceq_{\cF} N_1\leq N^+\ \&\  \|N_1\|\leq LS(\bT)\imply
q\restriction N_1$ does not fork over $N_0]$, we can get it for $\|N_1\|<\mu$, but does
it hold for $N_1=N^+$?  A central point is
\begin{description}
\item[$(*)$] Does $K$ satisfies amalgamation?
\end{description}
Again it seems that ${\cal P}^-(n)$-systems are called for.
Now Grossberg Shelah have started in the mid eighties to write a
paper, which solves the problem but with two
drawbacks. It says: if $\bT\subseteq L_{\kappa^*, \omega}$ has
arbitrarly large models, is categorical in $\chi^{+n}$ (for
$n<\omega$), $\chi\geq LS(\bT)$, and $2^{\chi^{+n}}<2^{\chi^{+n+1}}$
for $n<\omega$, then
$\bT$ is categorical in every $\lambda'>\chi$. So we need the
set theoretic assumption 
$$
\left(\forall \alpha<(2^{LS(\bT)})^+\right)\left(\exists
\chi\right)\left[\beth_\alpha\leq \chi\ \&\ \chi^{+n}\leq \lambda\ \&\
\bigwedge\limits_n 2^{\chi^{+n}}< 2^{\chi^{+n+1}}\right].
$$
\medskip

\noindent
4) If $|\bT|<\kappa^*$ we can do better, as $\Op(EM(I, \Phi))=EM(\Op(I),
\Phi)$, will discuss elsewhere.
\medskip

\noindent
5) Elsewhere we shall adopt what is done here to abstract elementary
class ${\cal K}$ categorical in $\lambda\geq \beth_{(2^{LS({\cal K})})^+}$
such that ${\cal K}_{<\lambda}$ has amalgamation.

\bibliographystyle{lit-plain}
\bibliography{listb,lista,listx}
\vfill
\eject

\shlhetal

\end{document}